\documentclass[12pt]{article}
\usepackage{amsmath}
\usepackage{amsthm}
\usepackage{amssymb}
\usepackage{tikz}
\usepackage{amsfonts}
\usepackage{epsf}
\usepackage{graphicx}
\usepackage{epstopdf}
\usepackage{appendix}
\usepackage{graphics}
\usepackage{graphicx}

\newtheorem{lemma}{Lemma}

\newtheorem{theorem}{Theorem}

\newtheorem{corollary}{Corollary}

\theoremstyle{remark}

\newtheorem{remark}{Remark}
\theoremstyle{definition}

\DeclareMathOperator\Ai{{Ai}}
\DeclareMathOperator\Bi{{Bi}}

\DeclareMathOperator\im{{Im}}
\numberwithin{equation}{section}

\numberwithin{equation}{section}

\newcounter{comment}
\def\sign{\mathbf{sgn}}

\begin{document}

\title{On the connection problem for the second Painlev\'e equation with large initial data}

\author{
Wen-Gao Long$^{1,}$\thanks{Wen-Gao Long:longwg@hnust.edu.cn,~$^{1}$School of Mathematics and Computational Science, Hunan University of Science and Technology, Xiangtan 411201, PR China}\quad and\quad Zhao-Yun Zeng$^{2,}$\thanks{Zhao-Yun Zeng(Corresponding author): zengzhaoyun@jgsu.edu.cn,~
$^{2}$School of Mathematics and Physics, Jinggangshan University, Ji'an 343009, PR China}}

\date{}

\maketitle

%%%%%%%%%%%%%%%%%%%%%%%%%%%%%%%%%%%%%%%%%%%%%%%%%%%%%%%%%%%%%%%%%%
\noindent \textbf{Abstract}\\
We consider two special cases of the connection problem for the homogenous second Painlev\'e equation (PII)
 using the method of {\it uniform asymptotics} proposed by Bassom {\it et al.}. We give an asymptotic
 classification of the real solutions of PII on the negative (positive) real axis with respect to their
 initial data. As by product, a rigorous proof of a property associated with the nonlinear eigenvalue problem
 of PII on the real axis, recently revealed by Bender and Komijani, is given by deriving the asymptotic
 behavior of the Stokes multipliers.

\vspace{3mm}

\vspace{3mm}

%%%%%%%%%%%%%%%%%%%%%%%%%%%%%%%%%%%%%%%%%%%%%%%%%%%%%%%%%%%%%%%%%%%%%%%
\section{Introduction}

The Painlev\'{e} equations are, in general, irreducible in the sense that their solutions cannot be expressed
in terms of elementary functions or classical special functions. Therefore, the study on the asymptotic
 behavior of the Painlev\'{e} functions and their connection problem is a difficult topic in
  Painlev\'{e} theory. In this paper, we focus on the connection problem of the homogenous second Painlev\'{e}
   (PII) equation, {\it i.e.} the following equation:
\begin{equation}\label{PII}
\frac{d^2q}{dt^2}=2q^3+tq,
\end{equation}
whose solutions are well studied. Although the ``classic connection problem'' for PII concerning solutions of
 \eqref{PII} satisfying the boundary condition $q(t)\to 0$ as $t\to+\infty$ was given by Hastings and McLeod
 forty years ago \cite{HM-1980}, the PII equation \eqref{PII} remains a research topic of great current
 interest due to its extensive applications in mathematics and
mathematical physics. Recall that every Painlev\'{e} transcendent admits a Riemann-Hilbert characterization
 through the Stokes multipliers. In the case of the homogenous PII equation, the Stokes multipliers are
 subject to the constraints
\begin{equation}\label{Stokes-multiplier-restriction}
  s_{k+3}=-s_k,~k\in\mathbb{Z},~{\rm and}~s_1-s_2+s_3+s_1s_2s_3=0,s_{1}=\overline{s_{3}},s_{2}=\overline{s_{2}}.
\end{equation}
By collecting the results in the literature \cite{BI-2012,Fokas book,HM-1980,SA-1981}, we know that, when
 $t\to-\infty$, the asymptotic behavior of the real-valued solutions of \eqref{PII} are divided into three
 types according to the Stokes multiplier $s_1$:
\begin{enumerate}
\item [(N1)] If $|s_{1}|<1$, the solutions are oscillatory on the negative real axis, and satisfying
\begin{equation}\label{eq-PII-asym-N1}
q(t)=d(-t)^{-\frac{1}{4}}\cos\left(\frac{2}{3}(-t)^{\frac{3}{2}}-\frac{3d^2}{4}\ln(-t)+\phi\right)+\mathcal{O}\left((-t)^{-\frac{7}{10}}\right)\quad
\end{equation}
as $t\to-\infty$, where
\begin{equation}\label{def-d-phi-N1}
d^2=-\frac{1}{\pi}\ln(1-s_{1}s_{3})\quad \text{and}\quad \phi=-\frac{3}{2}d^{2}\ln{2}+\arg\Gamma\left(\frac{1}{2}i d^2\right)-\arg{s_{1}};
\end{equation}
\item [(N2)] If $|s_{1}|=1$, the solutions behave like
\begin{equation}\label{eq-PII-asym-N2}
q(t)=-\epsilon \sqrt{-\frac{t}{2}}+h(-t)^{-\frac{1}{4}}e^{-\frac{2\sqrt{2}}{3}(-t)^{3/2}}\left(1+\mathcal{O}\left((-t)^{-\frac{3}{2}}\right)\right)
\end{equation}
as $t\to-\infty$, where
\begin{equation}\label{def-b-N2}
h=-\frac{s_{2}}{2^{7/4}\sqrt{\pi}}, \quad s_{1}=i\epsilon;
\end{equation}

\item [(N3)] If $|s_{1}|>1$, the solutions have an infinite number of poles on the negative real axis and
\begin{equation}\label{eq-PII-asym-N3}
q(t)=\frac{\sqrt{-t}}{\sin\left(\frac{2}{3}(-t)^{3/2}+\frac{3}{2}\beta\ln(-t)+\varphi\right)+\mathcal{O}\left((-t)^{-1/5}\right)}
\end{equation}
as $t\to-\infty$, where
\begin{equation}\label{def-beta-varphi-N3}
\beta=\frac{1}{2\pi}\ln(s_{1}s_{3}-1)-1, \quad \varphi=3\beta\ln{2}-\arg\Gamma\left(\frac{1}{2}+i\beta\right)-\arg{s_{1}}.
\end{equation}
\end{enumerate}
When $t\to+\infty$, the asymptotic behavior of the real-valued solutions of PII equation can be separated into two types according to the Stokes multiplier $s_2$:
\begin{enumerate}
\item [(P1)] If $s_{2}\neq 0$, the solutions (a two-parameter family) satisfy
\begin{equation}\label{eq-PII-asym-P1}
q(t)=\sigma\sqrt{\frac{t}{2}}\cot\left(\frac{\sqrt{2}}{3}t^{3/2}+\frac{3\gamma}{4}\ln{t}+\chi\right)+\mathcal{O}(t^{-1})
\end{equation}
as $t\to+\infty$ with $\sigma=\text{sgn}(s_{2})=\pm 1$ and
\begin{equation}\label{def-gamma-chi-P1}
 \gamma=\frac{1}{\pi}\ln|s_{2}| ,\quad \chi=\frac{7\gamma}{4}\ln{2}-\frac{1}{2}\arg\Gamma\left(\frac{1}{2}+i\gamma\right)-\frac{1}{2}\arg(1+s_{2}s_{3})+\frac{\pi}{2};
\end{equation}
\item [(P2)] If $s_{2}=0$, the solutions (a one-parameter family) behave like the Airy functions to leading order, {\it i.e.}
\begin{equation}\label{eq-PII-asym-P2}
q(t)=\kappa\Ai(t)(1+\mathcal{O}(t^{-3/2}))
\end{equation}
as $t\to+\infty$, where
\begin{equation}\label{def-kappa-P2}
\kappa=-\im s_{1}.
\end{equation}
\end{enumerate}

According to \eqref{eq-PII-asym-N1}-\eqref{def-kappa-P2}, one can build the connection formulas between the parameters in the asymptotic behaviors as $t\to+\infty$ and the ones as $t\to-\infty$. These connection formulas are well studied through the isomonodromy method \cite{Its-Kapaev-1988, Sule-1987}, the Deift-Zhou nonlinear steepest descent method \cite{BI-2012, DZ-1995} and the uniform asymptotics method \cite{APC}. See \cite{Fokas book} for a detailed summary of the connection formulas of PII and their applications.

The initial value problem of the inhomogenous PII equation with $(q(0),q'(0))=(0,0)$ was analyzed through the
 isomonodromy method in \cite{Kitaev-95}. To the best of our knowledge, it is the only case that the large $t$
  asymptotic behaviors are explicitly connected with the initial data. Even so, some heuristic numerical
 results of the PII equation should be noted. For instance, it can be seen from
    \cite[Figure 5]{Fornberg-Weideman-2014} that the real solutions of the PII equation have
    $n (n=1,2,\cdots)$ poles on the positive real axis when their initial data locate on certain curves marked $n^{+}$. This means that these solutions belong to type (P2).  Similarly, there exists a sequence of
    regions marked $n^{-}$ such that the solutions, whose initial data lies in these regions, possess $n$
    poles on the negative real axis. Then they are type (N1) solutions. Nevertheless, it is to be further
    investigated on how to theoretically describe these curves and regions on the $(q(0),q'(0))$ plane.

\begin{figure}
\setlength{\abovecaptionskip}{-0.5cm}
\setlength{\belowcaptionskip}{-0cm}
  % Requires \usepackage{graphicx}
  \includegraphics[width=17cm]{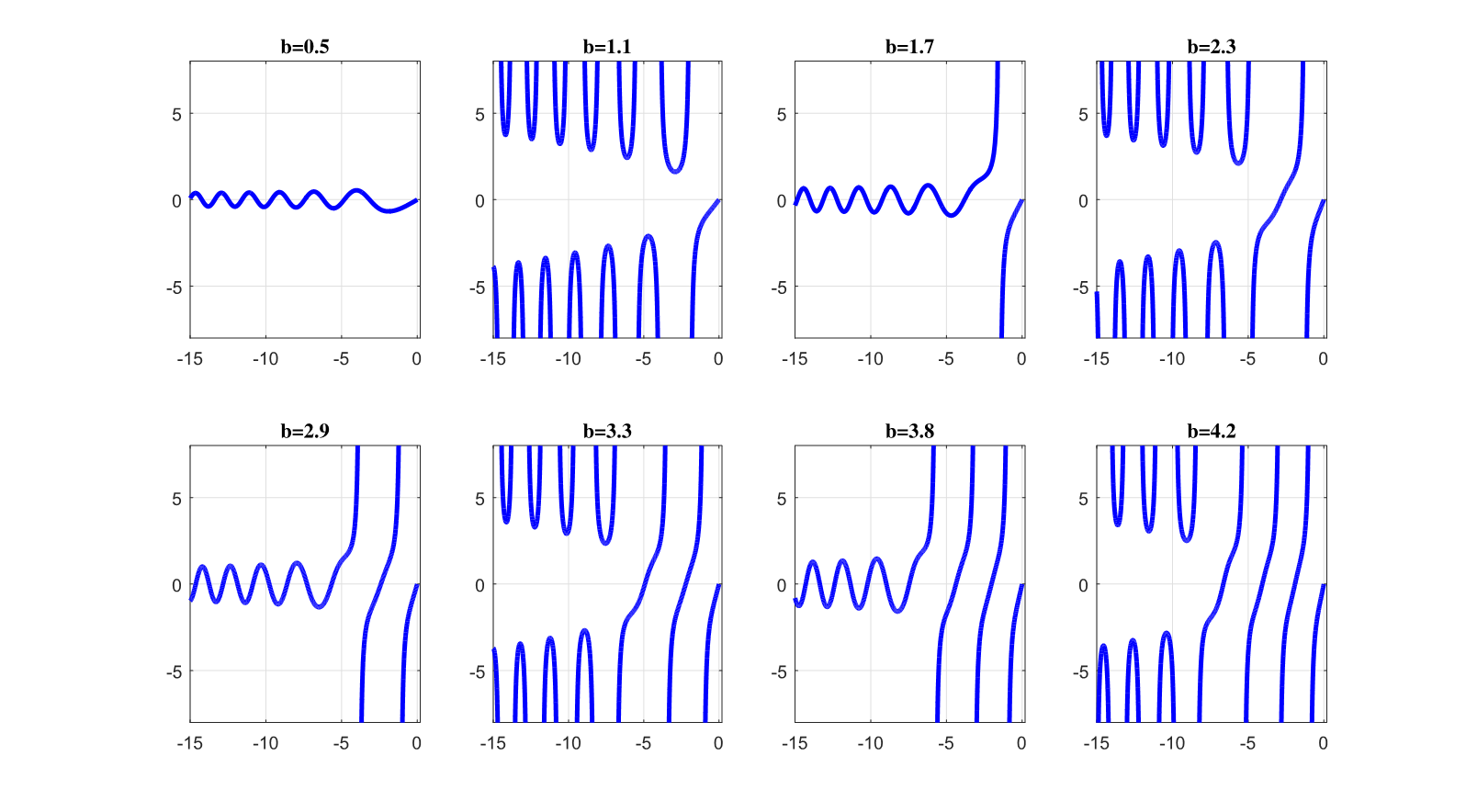}\\
  \caption{The type (N1) and (N3) solutions of PII with $q(0)=0$ and $q'(0)=b$ varies}\label{PII3types}
\end{figure}

Recently, the initial value problem of the homogenous PII equation with real initial data $(q(0),q'(0))$ have been numerically and analytically investigated by Bender {\it et al.}\cite{Bender-Komijani-2015, Bender-Komijani-Wang-2019}. They conclude that, when $q(0)=0$, there exist initial slopes $q'(0)=b_{n}$ that give rise to the solutions of type (N2), where
\begin{equation}\label{asymptotic-bn-En}
b_{n}\sim \left[\frac{3\sqrt{2\pi}\Gamma(\frac{3}{4})n}{\Gamma(\frac{1}{4})}\right]^{\frac{2}{3}} \quad\text{as}\quad n\rightarrow\infty.
\end{equation}
Moreover, a very interesting phenomenon appears as the initial slope varies. If $b_{2n}<q'(0)<b_{2n+1}$, then the solution passes through $n$ simple poles and then oscillates stably on the negative real axis (type (N1) solutions), while for $b_{2n-1}<q'(0)<b_{2n}$, the solution has infinite number of poles (type (N3) solutions); see Figure \ref{PII3types}. It is similar when $q'(0)=0$ and $q(0)$ varies; see \cite[Figures 14 and 15]{Bender-Komijani-2015}.
We note that parts of the argument in \cite{Bender-Komijani-2015, Bender-Komijani-Wang-2019} are based on numerical simulations (see \cite[p.\,8]{Bender-Komijani-2015}). In order to give a rigorous proof of their results and to build the connection formulas between the initial data and the large negative (positive) $t$ asymptotics, further analysis is needed.

Therefore, it is natural to consider the following initial value problem of PII
\begin{equation}\label{Initial-problem}
\begin{cases}
\frac{d^{2}q}{dt^2}=2q^3+tq, \\
  q(0)=a,~~q'(0)=b,
  \end{cases}
\end{equation}
where the initial data $a$ and $b$ are two real constants. The main task is to classify the PII solutions with respect to the initial data, or to build the connection formulas between the parameters involved in the asymptotic behaviors shown above in \eqref{eq-PII-asym-N1}--\eqref{def-kappa-P2} and the initial data $a$ and $b$. Actually, as is well known, it is an open problem proposed by Clarkson  on several occasions \cite{CPA2003, CPA2006, CPA2019}.

It seems difficult to completely solve Clarkson's open problem for PII, nevertheless it is possible to give an asymptotic classification of PII solutions with respect to the initial data. Similar analysis for the first Painlev\'{e} equation (PI) has been investigated by Long {\it et al.} in \cite{LongLLZ}. By fixing the initial value $q(0)$, they have classified the real PI solutions asymptotically with respect to the initial slope $b$. The main idea is to compute the asymptotic behavior of the Stokes multipliers when $b$ is large. Similarly, when fix the initial slope $b$, the real PI solutions are classified asymptotically with respect to the initial value $a$.

The main objective of this paper is to give an asymptotic classification of the real-valued solutions of the homogenous PII equation \eqref{PII} on the real axis with respect to the initial data. The key difference between the present paper and \cite{Bender-Komijani-2015, Bender-Komijani-Wang-2019,LongLLZ} is that we investigate how the asymptotic behaviors of the real-valued solutions of PII alternate when the initial data $a$ and $b$ are both varied. Note that $H(t)=\frac{1}{2}p^2(t)-\frac{1}{2}q^4(t)-\frac{1}{2}tq(t)^2$ is the Hamiltonian for the following Hamiltonian system of the PII equation
\begin{equation}
\left\{\begin{matrix}
\frac{dq}{dt}&=\frac{\partial H}{\partial p},\\
\frac{dp}{dt}&=-\frac{\partial H}{\partial q},
\end{matrix}\right.
\end{equation}
where $p=\frac{dq}{dt}$; see \cite[Table 4.2]{Kawai-Takei-2005-book}.
We regard $-2H(0)=a^4-b^2$ as a whole term and only assume that $a^4-b^2$ is large, instead of fixing $a$ or $b$.
Set $a=\xi^{1/3}A(\xi)$ and $b=\xi^{\frac{2}{3}}B(\xi)$, where $\xi$ is a large positive real parameter and $A(\xi), B(\xi)$ are both real and bounded. Since $a^4-b^2$ can be large positive or large negative, we shall consider the initial value problem \eqref{Initial-problem} in the following two cases:
\begin{equation}\label{initial-condition-two-cases}
\begin{aligned}
({\text i})&\quad A(\xi)^4-B(\xi)^2=1; \\
({\text {ii}})&\quad A(\xi)^4-B(\xi)^2=-1.
\end{aligned}
\end{equation}
In case (i), we obtain an asymptotic classification of the PII solutions on the positive real axis with respect to $a^4-b^2$, and in case (ii), we classify the PII solutions on the negative real axis similarly.
\begin{remark}
The shapes of the curves and regions in \cite[Figure 5]{Fornberg-Weideman-2014} also motivate us to regard
$a^4-b^2$ as a whole term. The scalings of $a$ and $b$ are prepared for the presentation of the main results
and the later analysis. Observing that $\xi^{\frac{4}{3}}=|a^4-b^2|$, one may find that the two cases in \eqref{initial-condition-two-cases} represent the two cases when $a^4-b^2\to+\infty$ and $-\infty$ respectively.
\end{remark}

The main technique in this paper is based on the method of {\it uniform asymptotics} introduced by Bassom {\it et al.} \cite{APC}, and further applied by \cite{LongLLZ, LongZZ, Qin1, WZ1, WZ2, Zeng, Zeng2}. {\color{blue}To apply this method,} we first briefly outline some important information of the monodromy theory of the PII equation. The reader is referred to \cite{Fokas book, MJ} for more details.
One of the Lax pairs for \eqref{PII} is the system of linear ordinary equations
\begin{equation}\label{PII:lax pair}
  \left\{\begin{array}{ll}
    \frac{\partial\Psi}{\partial\lambda}=\{-i(4\lambda^2+t+2q^2)\sigma_3+4\lambda q\sigma_1-2q'\sigma_2\}\Psi\\
 \frac{\partial\Psi}{\partial t}=(-i\lambda\sigma_3+q\sigma_1)\Psi.
  \end{array}\right.
\end{equation}
Here and below the prime denotes the derivative with respect to $t$, and
\begin{equation*}
  \sigma_1=\left(\begin{array}{cc}
 0&1\\
 1&0
 \end{array}\right),~~
 \sigma_2=\left(\begin{array}{cc}
 0&-i\\
 i&0
 \end{array}\right),~~
 \sigma_3=\left(\begin{array}{cc}
 1&0\\
 0&-1
 \end{array}\right)
\end{equation*}
are the standard Pauli spin matrices. %and $r(t)=\frac{dq}{dt}$.
A direct calculation shows that the compatibility condition $\Psi_{\lambda t}=\Psi_{t\lambda}$ implies \eqref{PII}.
There exist canonical solutions $\Psi_{k}(\lambda,t)$, defined in a neighborhood of the irregular singular point $\lambda=\infty$, of the first equation of the Lax pair, with the following asymptotics behavior
\begin{equation} \label{Psi boundary at infty1}
  \Psi_{k}(\lambda,t)=\left(I+O\left(\frac{1}{\lambda}\right)\right)\exp\left\{-\left(i\frac{4}{3}\lambda^3+it\lambda\right)\sigma_{3}\right\},
  ~\lambda\rightarrow\infty,~\lambda\in\Omega_k,
 \end{equation}
 where the canonical sectors are
 \begin{equation}\label{domain Omegak}
  \Omega_k=\left\{\lambda\in\mathbb{C}\Big| \frac{(k-2)\pi}{3}<\arg\lambda<\frac{k\pi}{3}\right\},~k\in\mathbb{Z}.
 \end{equation}
 These canonical solutions are related by
 \begin{equation}\label{Psi relation}
   \Psi_{k+1}(\lambda,x)=\Psi_k(\lambda,x)S_k,~\lambda\in\Omega_k\cap\Omega_{k+1},~k\in\mathbb{Z},
 \end{equation}
where $S_k$ are Stokes matrix defined by
\begin{equation}\label{Stokes matr}
  S_{2j}=\left(
                 \begin{array}{cc}
                   1 & s_{2j} \\
                   0 & 1 \\
                 \end{array}
               \right),\quad
  S_{2j-1}=\left(
                              \begin{array}{cc}
                                1 & 0\\
                                s_{2j-1} & 1 \\
                              \end{array}
                            \right), \quad j\in\mathbb{Z}.
\end{equation}
The Stokes multipliers $s_k,k\in\mathbb{Z}$ are independent of $\lambda$ and $t$ according to the isomonodromy condition.
As mentioned before, the Stokes multipliers are subject to the constraints \eqref{Stokes-multiplier-restriction}.

The rest of this paper is arranged as follows. In Section 2, we state our main results in two theorems. In Section 3, the main theorems are proved with two lemmas to be shown later. In Section 4, we apply the method of {\it uniform asymptotics} to calculate the Stokes multipliers when $t=0$ and $\xi\to+\infty$ in the two cases in \eqref{initial-condition-two-cases}, and then prove the lemmas in Section 3.  Some detailed calculations are left in the Appendixes.

\section{Main results}

When $A(\xi)^4-B(\xi)^2 = 1$ {\it i.e.} $a^{4}-b^{2}=\xi^{\frac{4}{3}}\to+\infty$, we have the following result.
\begin{theorem}\label{thm-case-I}
Let $q(t; a,b)$ be any {\color{blue}real-valued} on $\mathbb{R}$ solution of \eqref{Initial-problem}. Then there exists a sequence of curves
$$\Gamma_{n}:\quad a^4-b^2=f_{n}(a,b),\quad n=1,2,\cdots$$
on the $(a,b)$ plane such that $q(t; a,b)$ belongs to type (P2) solutions of PII when $(a,b)$ lies on these curves, where \begin{equation}\label{limit-connection-case-I-kappa}
\kappa:=\kappa_{n}=\frac{(-1)^{n}e^{2\xi E_{1}-2E_{2}+o(1)}}{A(\xi)}
\end{equation}
as $\xi\to+\infty$. As $n\to\infty$, we have
\begin{equation}\label{asym-fn}
f_{n}(a,b)=\left[\frac{n\pi-\frac{\pi}{2}+2E_{3}+o(1)}{2E_{1}}\right]^{\frac{4}{3}}.
\end{equation}
Moreover if $f_{n}(a,b)<a^4-b^2<f_{n+1}(a,b)$, then $q(t; a,b)$ belongs to type (P1) with $\sigma=\sign((-1)^{n+1} a)$ and
\begin{equation}\label{limit-connection-I-1}
\begin{split}
 \gamma&=\frac{1}{\pi}\left[-2\xi E_{1}+2E_{2}+\ln\left|2A(\xi)\cos\left[2\xi E_{1}-2E_{3}+o(1)\right]\right|+o(1)\right], \\ \chi&=\frac{7\gamma}{4}\ln{2}-\frac{1}{2}\arg\Gamma\left(\frac{1}{2}+i\gamma\right)-2\xi E_{1}+2E_{3}+o(1)
\end{split}\end{equation}
as $\xi\rightarrow+\infty$.
In the above formulas \eqref{asym-fn} and \eqref{limit-connection-I-1}, $E_{1},E_{2}$ and $E_{3}$ are defined by
\begin{equation}\label{def-E1-E2-E3}
\begin{split}
E_{1}&=\frac{1}{6}B\left(\frac{1}{2},\frac{1}{4}\right),\\
E_{2}&=\frac{1}{16}\frac{B(\xi)}{A(\xi)}B\left(\frac{1}{4},\frac{1}{2}\right)+\Re\mathcal{F}(\xi),\\
E_{3}&=\frac{1}{16}\frac{B(\xi)}{A(\xi)}B\left(\frac{1}{4},\frac{1}{2}\right)+\frac{\pi}{8}+\Im\mathcal{F}(\xi),
\end{split}\end{equation}
where $B(\cdot,\cdot)$ is the beta function and
$$\mathcal{F}(\xi):=\frac{2 A(\xi)^2-\left(B(\xi)/A(\xi)\right)^2}{8}\int_{\frac{1-i}{2}}^{\infty}\frac{1}{(s+\frac{B(\xi)}{2 A(\xi)})\sqrt{s^4+\frac 1 4}}ds.$$
Here and after, the branches of $\sqrt{s^4+1/4}$ are chosen such that $\arg(s-\alpha_{i})\in(-\pi,\pi)$, $i=1,2,3,4$, where $\alpha_1=\frac{1}{\sqrt{2}}e^{-\frac{\pi i}{4}}$, $\alpha_2=\frac{1}{\sqrt{2}}e^{\frac{\pi i}{4}}$, $\alpha_3=\frac{1}{\sqrt{2}}e^{\frac{3\pi i}{4}}$ and $\alpha_4=\frac{1}{\sqrt{2}}e^{-\frac{3\pi i}{4}}$.
\end{theorem}

Specially, in the case that the initial slope $y'(0)=b$ is fixed, noting that $\xi$ is a large parameter and $A(\xi)^4-B(\xi)^2=1$, we find that $B(\xi)\to 0$ and $A(\xi)\to \pm 1$. Moreover, a simple calculation yields $\mathcal{F}(\xi)=\frac{\pi}{8}+o(1)$ when $B(\xi)\to 0$. Hence, the following corollary is a direct consequence of Theorem \ref{thm-case-I} and is consistent with the corresponding result mentioned in \cite[p11]{Bender-Komijani-Wang-2019}.
\begin{corollary}
For any fixed $b\in\mathbb{R}$, there exists $0<a_{1}<a_{2}<\cdots<a_{n}<\cdots$ such that $q(t; a_{n},b)$ belong to type (P2) solutions of PII, and
\begin{equation}\label{a-n}
a_{n}^4-b^2=\left[\frac{n\pi+o(1)}{2E_{1}}\right]^{\frac{4}{3}}\quad \text{as} \quad  n\to\infty.
\end{equation}
\end{corollary}

\begin{figure}[h]
\centering
  % Requires \usepackage{graphicx}
  \includegraphics[width=12cm,height=12cm]{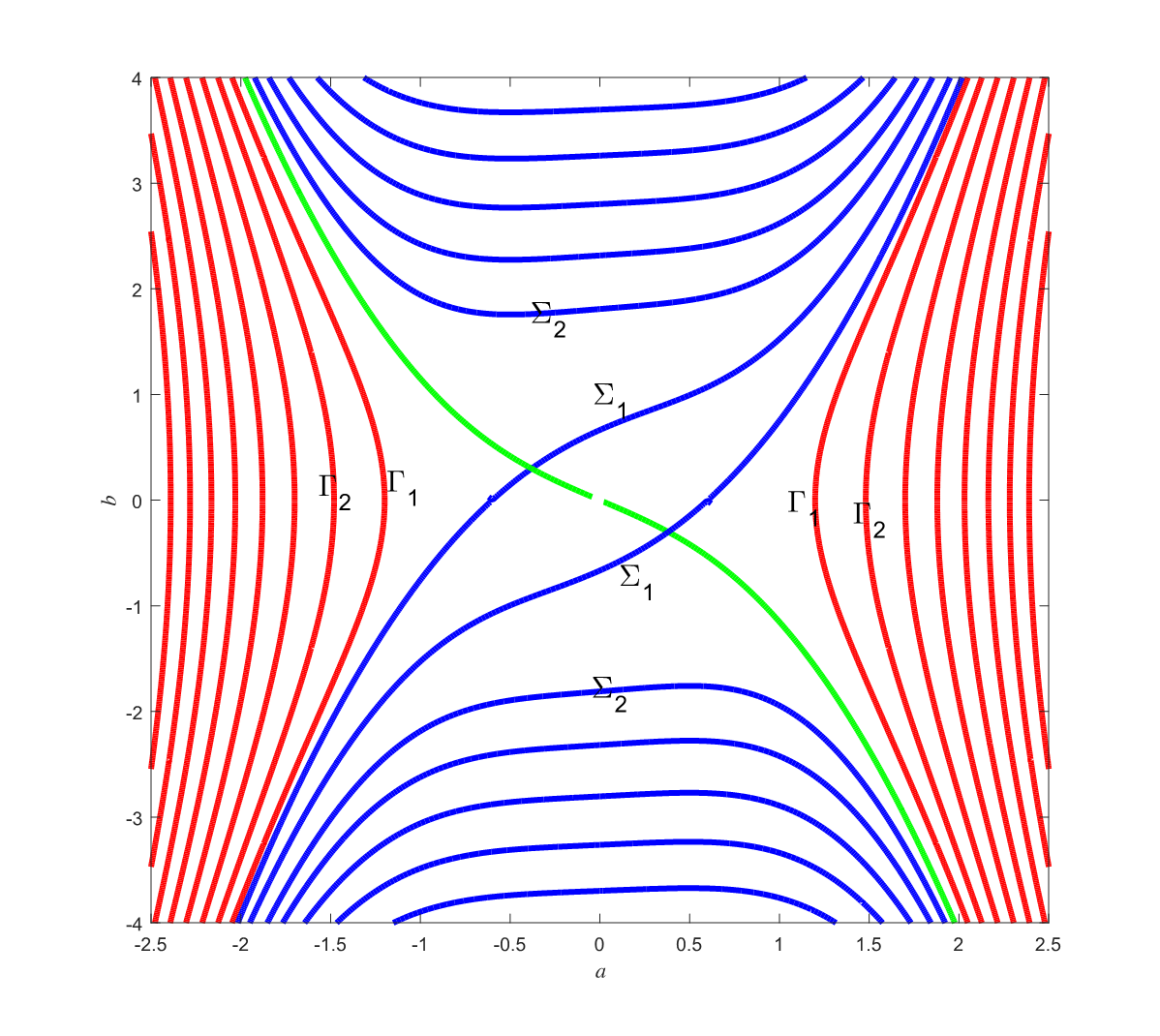}\\
  \caption{Sketch map of the curves $\Gamma_{n}$ and $\Sigma_{n}$; see the first subplot in \cite[Fig.5]{Fornberg-Weideman-2014}}\label{curves-gamma-sigma}
\end{figure}

\begin{remark}
According to \cite[Figure 5]{Fornberg-Weideman-2014}, there is a curve (the green curve in Figure \ref{curves-gamma-sigma}) such that the PII solutions, whose initial data locate on it, have no poles on the positive real axis. It means that these solutions belong to type (P2). Hence, one may regard this curve as $\Gamma_{0}$, we ignored it since the asymptotic form equation \eqref{asym-fn} is only valid for large $n$.
\end{remark}

When $A(\xi)^4-B(\xi)^2 =-1$ {\it i.e.} $a^{4}-b^{2}=-\xi^{\frac{4}{3}}\to-\infty$, similar results are obtained.
\begin{theorem}\label{thm-case-II}
Let $q(t; a,b)$ be any real-valued on $\mathbb{R}$ solution of \eqref{Initial-problem}. Then there exists a sequence of curves
$$\Sigma_{n}:\quad a^4-b^2=-g_{n}(a,b),\quad n=1,2,\cdots$$
on the $(a,b)$ plane such that $q(t; a,b)$ belongs to type (N2) solutions of PII equation when $(a,b)$ lies on these curves, where  $\epsilon=\sign((-1)^{n-1} b)$ and
\begin{equation}\label{h-n}
h:=h_{n}=-\frac{1}{2^{7/4}\sqrt{\pi}}B(\xi)e^{2\sqrt{2}E_{1}\xi-2F_{1}+o(1)}
\end{equation}
as $\xi\to+\infty$. As $n\to\infty$, we have
\begin{equation}\label{asym-gn}
g_{n}(a,b)=\left[\frac{n\pi-\pi+\sign(b)\cdot\frac{\pi}{2}-2F_{2}+o(1)}{2\sqrt{2}E_{1}}\right]^{\frac{4}{3}}.
\end{equation}
Moreover,
\begin{enumerate}
\item [(i)] if $g_{2n-1}(a,b)<-a^4+b^2<g_{2n}(a,b)$, then $q(t; a,b)$ belongs to type (N1) with
\begin{equation}\label{limit-connection-case-II-1}
\begin{split}
d^2&=-\frac{1}{\pi}\ln\left(\frac{2\cos(2\sqrt{2}\xi E_{1}+2F_{2})}{B(\xi)}\right)+\frac{2\sqrt{2}\xi E_{1}+2F_{1}}{\pi}+o(1),\\
\phi&=-\frac{3}{2}d^{2}\ln{2}+\arg\Gamma\left(\frac{1}{2}i d^2\right)+2\sqrt{2}\xi E_{1}+2F_{2}+o(1)
\end{split}
\end{equation}
as $\xi\to+\infty$.
\item [(ii)] if $g_{2n}(a,b)<-a^4+b^2<g_{2n+1}(a,b)$, then $q(t; a,b)$ belongs to type (N3), where
\begin{equation}\label{limit-connection-case-II-3}
\begin{split}
\beta&=\frac{1}{2\pi}\ln\left(-\frac{2\cos(2\sqrt{2}\xi E_{1}+2F_{2})}{B(\xi)}\right)-\frac{2\sqrt{2}\xi E_{1}+2F_{1}}{2\pi}-1, \quad \\
\varphi&=3\beta\ln{2}-\arg\Gamma\left(\frac{1}{2}+i\beta\right)+2\sqrt{2}\xi E_{1}+2F_{2}+o(1)
\end{split}
\end{equation}
as $\xi\to+\infty$.
\end{enumerate}
In the above formulas, $E_{1}$ is the same as in Theorem 1 and
\begin{equation}\label{def-F1-F2}
\begin{split}
F_{1}&=\frac{1}{2}\int_{\frac{1}{\sqrt{2}}}^{\infty}\frac{A(\xi)}{B(\xi)}\cdot \frac{1-A(\xi)^2}{2\frac{A(\xi)}{B(\xi)}s+1}\cdot \frac{1}{\sqrt{s^4-\frac{1}{4}}}ds-\frac{1}{2}\int_{\frac{1}{\sqrt{2}}}^{\infty}\frac{s\log(2\frac{A(\xi)}{B(\xi)}s+1)}{\left(s^2+\frac{1}{2}\right)\sqrt{s^4-\frac{1}{4}}}ds,\\
F_{2}&=\frac{i}{2}\int_{\frac{i}{\sqrt{2}}}^{\infty} \frac{A(\xi)}{B(\xi)}\cdot\frac{1+A(\xi)^2}{2\frac{A(\xi)}{B(\xi)}s+1}\cdot\frac{1}{\sqrt{s^4-\frac{1}{4}}}ds-\frac{i}{2}\int_{\frac{i}{\sqrt{2}}}^{\infty}\frac{s\log(2\frac{A(\xi)}{B(\xi)}s+1)}{\left(s^2-\frac{1}{2}\right)\sqrt{s^4-\frac{1}{4}}}ds+\frac{i\log{B(\xi)}}{2}.
\end{split}\end{equation}
Here and after, the branches of $\sqrt{s^4-1/4}$ are chosen such that $\arg(s-\hat{\alpha}_{i})\in(-\pi,\pi)$, $i=1,2,3,4$, where $\hat{\alpha}_1=\frac{1}{\sqrt{2}}$, $\hat{\alpha}_2=\frac{1}{\sqrt{2}}e^{\frac{\pi i}{2}}$, $\hat{\alpha}_3=-\frac{1}{\sqrt{2}}$ and $\hat{\alpha}_4=\frac{1}{\sqrt{2}}e^{-\frac{\pi i}{2}}$.
\end{theorem}

The following corollary is a direct consequence of Theorem \ref{thm-case-II} and is consistent with the corresponding result mentioned in \cite[p11]{Bender-Komijani-Wang-2019}.
\begin{corollary}
For any fixed $a\in\mathbb{R}$, there exist a sequence $0<b_{1}<b_{2}<\cdots<b_{n}<\cdots$ such that $q(t; a,b_{n})$ belong to type (N2) solutions of PII, and
\begin{equation}
a^4-b_{n}^2=-\left[\frac{n\pi-\frac{\pi}{2}+o(1)}{2\sqrt{2}E_{1}}\right]^{\frac{4}{3}}\quad \text{as} \quad n\to\infty.
\end{equation}
\end{corollary}

\begin{remark}
It should be noted that $E_{2}, E_{3}$ and $F_{1}, F_{2}$ are all functions of $a,b$. In fact, a combination of $\left|a^4-b^2\right|=\xi^{\frac{4}{3}}$ and $a=\xi^{\frac{1}{3}}A(\xi)$ yields $A(\xi)=a\left|a^4-b^2\right|^{-\frac{1}{4}}$. Similarly, one may find that $B(\xi)=b\left|a^4-b^2\right|^{-\frac{1}{2}}$. Moreover, note from \eqref{def-F1-F2} that $F_{2}(a,b)=F_{2}(-a,-b)+\pi$ for all $b>0$. Hence from \eqref{asym-gn}, we have $g_{n}(a,b)=g_{n}(-a,-b)$.

\end{remark}

\section{Proof of the main results}

According to the isomonodromy theory, a general idea to solve connection problems is to calculate the Stokes
 multipliers of a specific solution in the two specific situations to be connected. In the initial problem of
  PII,  this  means  to calculate all $s_{k}$  when  $t\rightarrow\pm\infty$ and $t=0$. When
   $t\rightarrow\pm\infty$, as stated above in \eqref{def-d-phi-N1}, \eqref{def-b-N2}, \eqref{def-beta-varphi-N3}, \eqref{def-gamma-chi-P1} and \eqref{def-kappa-P2},
   the Stokes multipliers are well known. When $t=0$, it seems much more difficult to find the
   exact values of $s_{k}$. However, inspired  by the ideas in Sibuya \cite{Sibuya-1967},  we are able to
    obtain  their  asymptotic approximations  in the special cases considered here, as a step forward. These
    approximations, together with \eqref{def-d-phi-N1}, \eqref{def-b-N2}, \eqref{def-beta-varphi-N3}, \eqref{def-gamma-chi-P1} and \eqref{def-kappa-P2}, suffice to prove our theorems.

\begin{lemma}\label{stokes-multipliers-case-I}
When $A(\xi)^4-B(\xi)^2=1$, the asymptotic behaviors of the Stokes multipliers corresponding to the PII solution $q(t; a,b)$ are
\begin{equation}\label{asymptotic-stokes-case-I}
\begin{split}
s_{2}&=-s_{-1}=-2A(\xi)e^{-2\xi E_{1}+2E_{2}+o(1)}\cos\left[2\xi E_{1}-2E_{3}+o(1)\right],\\
s_{1}&=\overline{s_{3}}=-\overline{s_{0}}=\frac{e^{2\xi(1+i)E_{1}-2(E_{2}+i E_{3})+o(1)}}{A(\xi)}
\end{split}
\end{equation}
as $\xi\to+\infty$, where $E_{1},E_{2},E_{3}$ are given in Theorem \ref{thm-case-I}.
\end{lemma}

\begin{corollary}\label{cor-case-I}
There exists a sequence of curves $\Gamma_{n}: a^4-b^2=f_{n}(a,b)$
on the $(a,b)$ plane such that the Stokes multipliers corresponding to $q(t; a,b)$ satisfy
$$s_{2}=s_{-1}=0,\quad s_{1}+s_{3}=0.$$
Moreover, we have \eqref{asym-fn} as $n\to\infty$.
\end{corollary}

\noindent\textbf{Proof of Corollary \ref{cor-case-I} and Theorem \ref{thm-case-I}} Regard $s_{2}$ as a function of $\xi$. Then, from \eqref{asymptotic-stokes-case-I}, we can define $f_{n}(a,b)$ by the $n$th zero of $s_{2}$ on the positive real axis. Moreover, we have the following two facts: (1) $s_{2}=0$ when $a^4-b^2=f_{n}(a,b)$; (2) $\arg{s_{1}}=n\pi-\frac{\pi}{2}-\arg{A(\xi)}+o(1)$ as $\xi\to+\infty$. It implies \eqref{asym-fn} immediately by substituting these two facts into \eqref{asymptotic-stokes-case-I}.

When $s_{2}=0$, using the constraints of the Stokes multipliers in \eqref{Stokes-multiplier-restriction}, we know that  $s_{1}=\overline{s_{3}}$ and $s_{1}+s_{3}=0$ . Hence, we obtain $\im s_{1}=\sign((-1)^{n+1}a)|s_{1}|$, which implies \eqref{limit-connection-case-I-kappa} immediately.
 When $s_{2}\neq 0$, it is readily seen from Lemma \ref{stokes-multipliers-case-I} that $\sign(s_{2})=\sign((-1)^{n+1}a)$ if $f_{n}(a,b)<a^4-b^2<f_{n+1}(a,b)$.  Hence, substituting \eqref{asymptotic-stokes-case-I} into \eqref{def-gamma-chi-P1}, we get \eqref{limit-connection-I-1}. This finishes the proof of Theorem \ref{thm-case-I}.

\ \\
Similarly, according to the following lemma and corollary, one can easily get Theorem \ref{thm-case-II}.
\begin{lemma}\label{stokes-multipliers-case-II}
When $A(\xi)^4-B(\xi)^2=-1$, the asymptotic behaviors of the Stokes multipliers corresponding to the PII solution $q(t; a,b)$ are
\begin{equation}
\begin{split}
s_{2}&=-s_{-1}=B(\xi)e^{2\sqrt{2}E_{1}\xi-2F_{1}+o(1)},\\
s_{1}&=\overline{s_{3}}=-\overline{s_{0}}=e^{-2\sqrt{2}\xi i E_{1}-2iF_{2}+o(1)}
\end{split}
\end{equation}
as $\xi\to+\infty$, where $E_{1},F_{1},F_{2}$ are given in Theorem \ref{thm-case-II}.
\end{lemma}

As an immediate consequence, we have the following corollary.
\begin{corollary}\label{cor-case-II}
There exists a sequence of curves $\Sigma_{n}: a^4-b^2=-g_{n}(a,b)$
on the $(a,b)$ plane such that the Stokes multipliers corresponding to $q(t; a,b)$ satisfy
$$s_{1}=\overline{s_{3}}=e^{-\sign(b)\frac{\pi i}{2}-(n-1)\pi i}.$$
Moreover, we have $g_{n}(a,b)=\left[\frac{(n-1)\pi+\sign(b)\cdot\frac{\pi}{2}-2F_{2}+o(1)}{2\sqrt{2}E_{1}}\right]^{\frac{4}{3}}$ as $n\to\infty$.
\end{corollary}

We leave the proofs of Lemmas \ref{stokes-multipliers-case-I} and \ref{stokes-multipliers-case-II} to the next section.

\section{Uniform asymptotics and proofs of the lemmas }\label{Uniform infty}
Let $\Psi_{k}=((\Phi_k)_1,(\Phi_k)_2)^T, k\in\mathbb{Z}$ be independent solutions of \eqref{PII:lax pair}, and set
\begin{equation}\label{define phi}
  \phi=\Big(4\lambda q-2iq'\Big)^{-\frac{1}{2}}(\Phi_k)_2,
\end{equation}
we get from \eqref{PII:lax pair} the second-order linear differential equation for $\phi(\lambda)$
\begin{align}\label{Sec order equ}
  \frac{d^2\phi}{d\lambda^2}
  &=\left\{-(4\lambda^2+t+2q^2)^2+8i\lambda+(4\lambda q+2iq')(4\lambda q-2iq')\right.\nonumber\\
  &\left.\quad-\frac{4iq(4\lambda^2+t+2q^2)}{4\lambda q-2iq'}+
  \frac{3}{4}\frac{(4q)^2}{(4\lambda q-2iq')^2}\right\}\phi.
\end{align}
When $t=0$, equation \eqref{Sec order equ} is simplified to
\begin{align}\label{Sec order equ2}
  \frac{d^2\phi}{d\lambda^2}&=\Big\{-\left(16\lambda^4+4a^4-4b^2\right)+4i\lambda+\frac{2b}{a}+\frac{i((b/a)^2 - 2a^2 )}{\lambda-\frac{ib}{2a}}+\frac{3}{4}\frac{1}{(\lambda-\frac{ib}{2a})^2}
  \Big\}\phi.
\end{align}
Our task is to derive the uniform asymptotic approximation of the solutions of \eqref{Sec order equ2} as $\xi\rightarrow+\infty$. For convenience, we assume that $A(\xi)\cdot B(\xi)>0$, and the analysis is similar when $A(\xi)\cdot B(\xi)<0$ and only justified in a remark; (see Remark \ref{remark-case-I-negative}).

\subsection{Case one: $A(\xi)^4-B(\xi)^2=1$}
Make the scaling
\begin{equation}\label{scaling-caseI}
\lambda=-i\xi^{\frac{1}{3}}\eta
\end{equation}
and set $Y(\eta,\xi)=\phi(\lambda)$. A straightforward calculation from \eqref{Sec order equ2} gives
\begin{equation}
\begin{split}\label{Second order equation asy-I}
  \frac{d^2Y}{d\eta^2}&=\xi^2\left\{16(\eta^4+\frac14)-\frac{4\eta}{\xi}-\frac{2B(\xi)}{\xi A(\xi)}-\frac{1}{\xi}\cdot\frac{2 A(\xi)^2-\left(\frac{B(\xi)}{A(\xi)}\right)^2}{\eta +\frac{B(\xi)}{2A(\xi)}}+g(\eta,\xi)\right\}Y\\
  &:=\xi^2G(\xi,\eta)Y(\eta,\xi),
\end{split}
\end{equation}
where $g(\eta,\xi)=\mathcal{O}\left(\frac{1}{\xi^{2}}\right)$ as $\xi\rightarrow+\infty$ uniformly for all $\eta$ bounded away from $-\frac{B(\xi)}{2A(\xi)}$.
For large $\xi$, it follows from equation \eqref{Second order equation asy-I} that there are four simple turning
points, say $\eta_j$, $j=1,2,3,4$. These turning points are near $\alpha_1=\frac{1}{\sqrt{2}}e^{-\frac{\pi i}{4}}$, $\alpha_2=\frac{1}{\sqrt{2}}e^{\frac{\pi i}{4}}$, $\alpha_3=\frac{1}{\sqrt{2}}e^{\frac{3\pi i}{4}}$ and $\alpha_4=\frac{1}{\sqrt{2}}e^{-\frac{3\pi i}{4}}$ respectively with $|\eta_j-\alpha_j|=\mathcal{O}\left(\xi^{-1}\right),~j=1,2,3,4$.
Moreover, it also follows from \eqref{Second order equation asy-I} that $\frac{1}{\eta_j-\alpha_j}=\mathcal{O}\left(\xi\right),~j=1,2,3,4$ as $\xi\to+\infty$.

\begin{figure}[!h]
  % Requires \usepackage{graphicx}
  \centering
  \includegraphics[width=6cm]{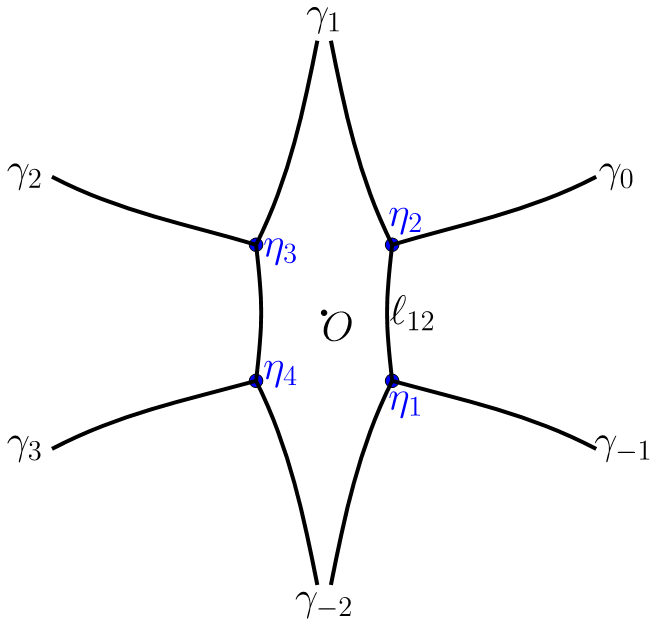}\\
  \caption{Stokes geometry of $G(\eta,\xi)d\eta^2$}\label{fig-stokes-I}
\end{figure}

According to \cite[Figure 5]{Eremenko}, the limiting state of the Stokes geometry of the quadratic differential $G(\eta,\xi)d\eta^2$ as $\xi\to+\infty$
 is described in Figure \ref{fig-stokes-I}. Therefore, following the main theorem in \cite{APC}, we can approximate the solutions of
 \eqref{Second order equation asy-I} in terms of the Airy functions.
Actually, if we define two conformal mappings $\vartheta(\eta)$ and $\zeta(\eta)$ by
\begin{equation}\label{Define vartheta}
\int_0^\vartheta s^{1/2}ds=\int_{\eta_1}^{\eta}G(s,\xi)^{1/2}d s,
\end{equation}
and
\begin{equation}\label{Define zeta}
  \int_0^{\zeta} s^{1/2}ds=\int_{\eta_2}^\eta G(s,\xi)^{1/2}ds,
\end{equation}
respectively in the neighborhoods of $\eta=\eta_1$ and $\eta=\eta_2$, then the following lemma is a direct consequence of \cite[Theorem 2]{APC}.
%%%%%%%%%%%%%%%%%%%%%%%%%%%%%%%%%%%%%%%%%%%%%%%%%%%%%%%%%%%%%%%%%%%%%%%%%%%%%%%%%%%%%%%%%%%%%%%%%%%%%%%%%%%%%%%%%%%%%%%%%%
\begin{lemma}\label{Uniform asymptotic theorem}
~Given any solution $Y(\eta,\xi)$ of \eqref{Second order equation asy-I},
there exist constants $C_1,~C_2$ and $\tilde{C}_1,~\tilde{C}_2$ such that,
 \begin{equation}\label{Uniform appr phi1}
   Y(\eta,\xi)=\Big(\frac{\vartheta}{G(\eta,\xi)}\Big)^{\frac{1}{4}}
\Big\{[C_1+o(1)]{\rm Ai}(\xi^{\frac{2}{3}}\vartheta)
+[C_2+o(1)]{\rm Bi}(\xi^{\frac{2}{3}}\vartheta)\Big\},~\xi\to +\infty
 \end{equation}
 uniformly for $\eta$ on any two adjacent Stokes lines emanating from $\eta_1$; and
\begin{equation}\label{Uniform appr phi2}
  Y(\eta,\xi)=\Big(\frac{\zeta}{G(\eta,\xi)}\Big)^{\frac{1}{4}}
\Big\{[\tilde{C}_1+o(1)]{\rm Ai}(\xi^{\frac{2}{3}}\zeta)
+[\tilde{C}_2+o(1)]{\rm Bi}(\xi^{\frac{2}{3}}\zeta)\Big\},~\xi\to +\infty
 \end{equation}
 uniformly for $\eta$ on any two adjacent Stokes lines emanating from $\eta_2$.
\end{lemma}
To calculate the connection matrices $S_{-1}$ and $S_{0}$ as $\xi\rightarrow+\infty$, it is necessary to clarify the asymptotic behavior of
 $\vartheta(\eta)$ and $\zeta(\eta)$ as $\xi,\eta\to\infty$.
%%%%%%%%%%%%%%%%%%%%%%%%%%%%%%%%%%%%%%%%%%%%%%%%%%%%%%%%%%%%%%%%%%%%
\begin{lemma}\label{lem relation zeta and eta infty}
Let $B(\cdot,\cdot)$ be the Euler's beta function. For large $\xi$ and $\eta$ with $\eta\gg\xi$, we have
\begin{equation}\label{eq-vartheta(eta)-infty}
  \frac{2}{3}\vartheta^{\frac{3}{2}}=\frac{4}{3} \eta^3+(1+i)E_{1}-\frac{1}{2\xi}\ln(2\eta)-\frac{E_{2}+iE_{3}}{\xi}+o\left(\frac{1}{\xi}\right),
\end{equation}
and
\begin{equation}\label{eq-zeta(eta)-infty}
  \frac{2}{3}\zeta^{\frac{3}{2}}=\frac{4}{3} \eta^3+(1-i)E_{1}-\frac{1}{2\xi}\ln(2\eta)-\frac{E_{2}-iE_{3}}{\xi}+o\left(\frac{1}{\xi}\right),
\end{equation}
where $E_{1}=\frac{1}{6}B(\frac12,\frac14)$, $E_{2}=\frac{B(\xi)}{16A(\xi)}B\left(\frac{1}{4},\frac{1}{2}\right)+\Re\mathcal{F}(\xi)$,  $E_{3}=\frac{B(\xi)}{16A(\xi)}B\left(\frac{1}{4},\frac{1}{2}\right)+\frac{\pi}{8}+\Im\mathcal{F}(\xi)$, and
$$\mathcal{F}(\xi):=\frac{2 A(\xi)^2-\left(B(\xi)/A(\xi)\right)^2}{8}\int_{\frac{1-i}{2}}^{\infty}\frac{1}{\left(s+\frac{B(\xi)}{2 A(\xi)}\right)\sqrt{s^4+\frac 1 4}}ds.$$
Specially, when $B(\xi)=0$, we find that $A(\xi)^2=1$. Then $E_{2}=0$ and $E_{3}=\frac{\pi}{4}$.
\end{lemma}
The proof of Lemma \ref{lem relation zeta and eta infty} will be given in Appendix A. Now we turn to the proof of Lemma \ref{stokes-multipliers-case-I}.

\textbf{Proof of Lemma \ref{stokes-multipliers-case-I}.}
According to \cite{APC}, in order to calculate the Stokes multipliers, one needs to know the uniform
asymptotic behaviors of $Y(\eta,\xi)$ on two adjacent Stokes lines. Therefore the uniform asymptotic behaviors
 of the Airy functions is required. In fact, according to \cite[Eqs.(9.2.12), (9.7.5), (9.2.10)]{OL}
 (see also \cite[Eqs.(4.5), (4.7)]{LongLLZ}), we know that
\begin{eqnarray}\label{eq-behavior-Ai}
\left\{\begin{aligned}
\Ai(z)\sim&\frac{1}{2\sqrt{\pi}}z^{-\frac{1}{4}}e^{-\frac{2}{3}z^{\frac{3}{2}}}, &\arg z\in \left(-\pi,\pi\right);\\
\Ai(z)\sim&\frac{1}{2\sqrt{\pi}}z^{-\frac{1}{4}}e^{-\frac{2}{3}z^{\frac{3}{2}}}+\frac{i}{2\sqrt{\pi}}z^{-\frac{1}{4}}e^{\frac{2}{3}z^{\frac{3}{2}}}, &\arg z\in\left(\frac{\pi}{3},\frac{5\pi}{3}\right);\\
\Ai(z)\sim&\frac{1}{2\sqrt{\pi}}z^{-\frac{1}{4}}e^{-\frac{2}{3}z^{\frac{3}{2}}}-\frac{i}{2\sqrt{\pi}}z^{-\frac{1}{4}}e^{\frac{2}{3}z^{\frac{3}{2}}}, &\arg z\in\left(\frac{-5\pi}{3},\frac{-\pi}{3}\right)
\end{aligned}\right.
\end{eqnarray}
and
\begin{equation}\label{eq-behavior-Bi}
\left\{\begin{aligned}
\Bi(z)\sim&\frac{1}{\sqrt{\pi}}z^{-\frac{1}{4}}e^{\frac{2}{3}z^{\frac{3}{2}}}+\frac{i}{2\sqrt{\pi}}z^{-\frac{1}{4}}e^{-\frac{2}{3}z^{\frac{3}{2}}},&\arg z \in \left(-\frac{\pi}{3},\pi\right);\\
\Bi(z)\sim&\frac{1}{\sqrt{\pi}}z^{-\frac{1}{4}}e^{\frac{2}{3}z^{\frac{3}{2}}}-\frac{i}{2\sqrt{\pi}}z^{-\frac{1}{4}}e^{-\frac{2}{3}z^{\frac{3}{2}}},&\arg z \in \left(-\pi,\frac{\pi}{3}\right);\\
\Bi(z)\sim&\frac{1}{2\sqrt{\pi}}z^{-\frac{1}{4}}e^{\frac{2}{3}z^{\frac{3}{2}}}+\frac{i}{2\sqrt{\pi}}z^{-\frac{1}{4}}e^{-\frac{2}{3}z^{\frac{3}{2}}},&\arg z\in\left(\frac{\pi}{3},\frac{5\pi}{3}\right)
\end{aligned}\right.
\end{equation}
as $z\rightarrow\infty$.

To derive $s_{-1}$, we should know the uniform asymptotic behaviors of $Y(\eta,\xi)$ on the two Stokes lines tending to infinity with $\arg\lambda\sim-\frac{2\pi}{3}$ and $\arg\lambda\sim-\frac{\pi}{3}$. In view of the transformation $\lambda=-i\xi^{\frac{1}{3}}\eta$, it is readily seen that, on these two lines, $\arg\eta\sim\mp\frac{\pi}{6}$ as $\eta\rightarrow\infty$. However, one may find in Figure \ref{fig-stokes-I} that the two adjacent Stokes lines $\gamma_{-1}$ and $\gamma_{0}$ with $\arg\eta\sim \pm \frac{\pi}{6}$ emanate from two different turning points. Hence we should build the relations between $(\Ai(\xi^{\frac{2}{3}}\vartheta),\Bi(\xi^{\frac{2}{3}}\vartheta))$ and $(\Ai(\xi^{\frac{2}{3}}\zeta),\Bi(\xi^{\frac{2}{3}}\zeta))$. Actually, it can be done by matching them on the Stokes line $\ell_{12}$ joining $\eta_{1}$ and $\eta_{2}$ in a similar way as in \cite[p518-p519]{LongLLZ}. Precisely, we have
\begin{equation}\label{eq-zeta-omega-relation}
\zeta^{\frac{1}{4}}\left(\Ai(\xi^{\frac{2}{3}}\zeta),\Bi(\xi^{\frac{2}{3}}\zeta)\right)\sim\vartheta^{\frac{1}{4}}\left(\Ai(\xi^{\frac{2}{3}}\vartheta),\Bi(\xi^{\frac{2}{3}}\vartheta)\right)
\left(\begin{matrix}
e^{Q(\xi)}&-i\left[e^{Q(\xi)}+e^{-Q(\xi)}\right]\\[0.2cm]
0&e^{-Q(\xi)}
\end{matrix}\right)
\end{equation}
as $\xi\rightarrow+\infty$, uniformly for $\eta\in \gamma_{0}\cup\gamma_{1}\cup \ell_{12}$, where
\begin{equation}\label{Q(xi)}
Q(\xi)=\xi\int_{\eta_{1}}^{\eta_{2}}G(s,\xi)^{\frac{1}{2}}ds
\end{equation}
and the branches are chosen such that $\arg(s-\eta_{1,2})\in(-\pi,\pi)$. Note that
\begin{equation}\label{Q(xi)-seperate}
\int_{\eta_{1}}^{\eta_{2}}G(s,\xi)^{\frac{1}{2}}ds=\int_{\eta_{1}}^{\eta}G(s,\xi)^{\frac{1}{2}}ds-\int_{\eta_{2}}^{\eta}G(s,\xi)^{\frac{1}{2}}ds.
\end{equation}
Then a combination of \eqref{Define vartheta}, \eqref{Define zeta}, \eqref{eq-vartheta(eta)-infty} and \eqref{eq-zeta(eta)-infty} yields
\begin{equation}\label{Q(xi)-E1-E3}
 Q(\xi)=i(2\xi E_{1}-2E_{3})+o(1)\quad \text{as}\quad \xi\rightarrow+\infty.
\end{equation}

Now we begin to calculate $s_{-1}$ via $\Psi_{0}(\lambda)=\Psi_{-1}(\lambda)S_{-1}$ and $s_{0}$ via $\Psi_{1}(\lambda)=\Psi_{0}(\lambda)S_{0}$.
If $\lambda\rightarrow\infty$ with $\arg\lambda\sim -\frac{\pi}{3}$, then $\arg\eta\sim\frac{\pi}{6}$, which
implies that $\arg\zeta\sim\frac{\pi}{3}$. Hence, substituting (\ref{eq-zeta(eta)-infty}) into
 \cite[Eq.(3.11) and Eq.(3.13)]{LongLLZ}, noting the transformation $\lambda=-i\xi^{\frac{1}{3}}\eta$ and
  observing the definition of $G(\eta,\xi)$ in \eqref{Second order equation asy-I}, we obtain
\begin{equation}\label{eq-Phi-pi/6}
\begin{cases}
\sqrt{2\lambda a-ib}\left(\frac{\zeta}{G(\eta,\xi)}\right)^{\frac{1}{4}}\Ai(\xi^{\frac{2}{3}}\zeta)&\sim c_{1}e^{\frac{4i}{3}\lambda^{3}},\\
\sqrt{2\lambda a-ib}\left(\frac{\zeta}{G(\eta,\xi)}\right)^{\frac{1}{4}}\Bi(\xi^{\frac{2}{3}}\zeta)&\sim ic_{1}e^{\frac{4i}{3}\lambda^{3}}+2c_{2}\frac{a}{2\lambda}e^{-\frac{4i}{3}\lambda^{3}},
\end{cases}
\end{equation}
where
\begin{equation}\label{eq-c1-c2-caseI}
c_{1}=\frac{a^{\frac{1}{2}}e^{-\frac{\pi i}{4}-\xi(1-i)E_{1}+(E_{2}-i E_{3})+o(1)}}{2\sqrt{\pi}},\qquad c_{2}= -\frac{\xi^{\frac{1}{3}}e^{\frac{\pi i}{4}+\xi(1-i)E_{1}-(E_{2}-i E_{3})+o(1)}}{2\sqrt{\pi}a^{\frac{1}{2}}}.
\end{equation}
as $\xi\rightarrow+\infty$. In view of the fact that $(\Phi_k)_{2}=((\Phi_k)_{21},(\Phi_k)_{22})$ and according to \cite[Eq.(8.1.7)]{Fokas book}, we know that
\begin{equation}\label{eq-asy-Phi21-Phi22}
(\Phi_0)_{21}\sim \frac{a}{2\lambda}e^{-\frac{4i}{3}\lambda^{3}},\qquad(\Phi_0)_{22}\sim e^{\frac{4i}{3}\lambda^{3}},\qquad \lambda\in\Omega_{k}
\end{equation}
as $\lambda\rightarrow\infty$ with $\arg \lambda=\arg(-i\eta)\sim-\frac{\pi}{3}$. Combining \eqref{eq-Phi-pi/6} and \eqref{eq-asy-Phi21-Phi22}, we have
\begin{equation}\label{Phi21-Phi22-to-Airy-gamma0}
\left((\Phi_0)_{21},(\Phi_0)_{22}\right)=\sqrt{2\lambda a-ib}\left(\frac{\zeta}{G(\eta,\xi)}\right)^{\frac{1}{4}}\left(\Ai(\xi^{\frac{2}{3}}\zeta),\Bi(\xi^{\frac{2}{3}}\zeta)\right)
\left(\begin{matrix}-\frac{i}{2c_{2}} & \frac{1}{c_{1}} \\ \frac{1}{2c_{2}}& 0\end{matrix}\right)
\end{equation}
as $\xi\rightarrow+\infty$. Here and after, the $c_{j}$'s in (\ref{Phi21-Phi22-to-Airy-gamma0}) are not equal but asymptotically equal to the corresponding ones in (\ref{eq-c1-c2-caseI}) as $\xi\rightarrow+\infty$. By a little abuse of notations, we use the same symbol for the $c_{j}$'s in these two formulas, since we only care about the asymptotic behavior of the Stokes multipliers.

If $\lambda\rightarrow\infty$ with $\arg\lambda\sim -\frac{2\pi}{3}$, then $\arg\eta\sim-\frac{\pi}{6}$, which implies that $\arg\vartheta\sim-\frac{\pi}{3}$. Using a similar argument in the derivation of \eqref{Phi21-Phi22-to-Airy-gamma0},
we get
\begin{equation}\label{Phi21-Phi22-to-Airy-gamma1}
((\Phi_{-1})_{21},(\Phi_{-1})_{22})=\sqrt{2(\lambda-a)}\left(\frac{\theta}{G(\eta,\xi)}\right)^{\frac{1}{4}}(\Ai(\xi^{\frac{2}{3}}\theta),\Bi(\xi^{\frac{2}{3}}\omega))
\left(\begin{matrix}
\frac{i}{2\tilde{c}_{2}}&\frac{1}{\tilde{c}_{1}}\\[0.2cm]
\frac{1}{2\tilde{c}_{2}}&0
\end{matrix}\right),
\end{equation}
where
\begin{equation}\label{eq-tilde-c1-c2}
\tilde{c}_{1}=\frac{a^{\frac{1}{2}}e^{-\frac{\pi i}{4}-\xi(1+i)E_{1}+(E_{2}+i E_{3})}}{2\sqrt{\pi}},\qquad
\tilde{c}_{2}=-\frac{\xi^{\frac{1}{3}}e^{\frac{\pi i}{4}+\xi(1+i)E_{1}-(E_{2}+i E_{3})}}{2\sqrt{\pi}a^{\frac{1}{2}}}.
\end{equation}
as $\xi\rightarrow+\infty$. A combination of \eqref{eq-zeta-omega-relation}, \eqref{Phi21-Phi22-to-Airy-gamma0}, \eqref{Phi21-Phi22-to-Airy-gamma1} and the fact $(\Phi_{0})_{2}=(\Phi_{-1})_{2}S_{-1}$ yields
\begin{equation}\label{eq-Stokes-matrix-S-1}
\begin{split}
S_{-1}&=
\left(\begin{matrix}
\frac{i}{2\tilde{c}_{2}}&\frac{1}{\tilde{c}_{1}}\\[0.2cm]
\frac{1}{2\tilde{c}_{2}}&0
\end{matrix}\right)^{-1}
\left(\begin{matrix}
e^{Q(\xi)}&-i\left[e^{Q(\xi)}+e^{-Q(\xi)}\right]\\[0.2cm]
0&e^{-Q(\xi)}
\end{matrix}\right)
\left(\begin{matrix}
-\frac{i}{2c_{2}} & \frac{1}{c_{1}} \\
 \frac{1}{2c_{2}}& 0
 \end{matrix}\right)\\
& =\left(\begin{matrix}1&0\\
 -2i\frac{\tilde{c}_{1}}{c_{2}}\cos\left[-iQ(\xi)\right]&0\end{matrix}\right).
 \end{split}
\end{equation}
To get the last equality, one may note the fact that $\frac{\tilde{c}_{2}}{c_{2}}=e^{Q(\xi)}$ which can be directly derived from \eqref{eq-c1-c2-caseI}, \eqref{eq-tilde-c1-c2} and \eqref{Q(xi)-E1-E3}. Finally, substituting the explicit expressions of $\tilde{c}_{1}$ and $c_{2}$ into \eqref{eq-Stokes-matrix-S-1}, we obtain $s_{-1}=2A(\xi)e^{-2\xi E_{1}+2E_{2}+o(1)}\cos\left[2\xi E_{1}-2E_{3}+o(1)\right]$ as $\xi\rightarrow+\infty$.

For the derivation of $s_{0}$, we also need the uniform asymptotic behavior of $Y(\eta,\xi)$ on the Stokes line $\gamma_{1}$. In this case, $\arg\lambda=\arg(-i\eta)\sim 0$ as $|\eta|\rightarrow+\infty$. In view of \eqref{eq-zeta(eta)-infty}, we find that $\arg\zeta\sim \pi$.
Using the asymptotic behaviors of the Airy functions $\Ai(z),\Bi(z)$ in \eqref{eq-behavior-Ai} and \eqref{eq-behavior-Bi}, it is readily seen that
\begin{equation}
\begin{cases}
\sqrt{2\lambda a-ib}\left(\frac{\zeta}{G(\eta,\xi)}\right)^{\frac{1}{4}}\Ai(\xi^{\frac{2}{3}}\zeta)&\sim c_{1}e^{\frac{4i}{3}\lambda^{3}}+ic_{2}\frac{ia}{2\lambda}e^{-\frac{4i}{3}\lambda^{3}},\\
\sqrt{2\lambda a-ib}\left(\frac{\zeta}{G(\eta,\xi)}\right)^{\frac{1}{4}}\Bi(\xi^{\frac{2}{3}}\zeta)&\sim ic_{1}e^{\frac{4i}{3}\lambda^{3}}+c_{2}\frac{ia}{2\lambda}e^{-\frac{4i}{3}\lambda^{3}}.
\end{cases}
\end{equation}
In a similar way of the above argument for deriving \eqref{Phi21-Phi22-to-Airy-gamma0} or \eqref{Phi21-Phi22-to-Airy-gamma1}, we obtain
\begin{equation}\label{eq-Phi-3pi/5-1}
\begin{aligned}
((\Phi_{1})_{21},(\Phi_{1})_{22})=\sqrt{2\lambda a-ib}\left(\frac{\zeta}{G(\eta,\xi)}\right)^{\frac{1}{4}}(\Ai(\xi^{\frac{2}{3}}\zeta),\Bi(\xi^{\frac{2}{3}}\zeta))
\left(\begin{matrix}                                                                                                                                                                                          -\frac{i}{2c_{2}}&\frac{1}{2c_{1}}\\[0.2cm]                                                                                                                                                                                 \frac{1}{2c_{2}}&-\frac{i}{2c_{1}}                                                                                                                                                                       \end{matrix}\right).
\end{aligned}
\end{equation}
Matching this result to the asymptotic behaviors of $(\Phi_{1})_{21},(\Phi_{1})_{22}$ in \eqref{eq-asy-Phi21-Phi22}, one can immediately get
\begin{equation}\label{stokes-matrix-s1}
((\Phi_{1})_{21},(\Phi_{1})_{22})=((\Phi_{0})_{21},(\Phi_{0})_{22})\left(\begin{matrix}
                                                             1 & -i\frac{c_{2}}{c_{1}} \\
                                                             0 & 1
                                                           \end{matrix}\right).
\end{equation}
where $c_{1}$ and $c_{2}$ are given in \eqref{eq-c1-c2-caseI}.  This implies $s_{0}=-i\frac{c_{2}}{c_{1}}=-\frac{e^{2\xi(1-i)E_{1}-2(E_{2}-i E_{3})+o(1)}}{A(\xi)}$ as $\xi\rightarrow+\infty$.

Finally, in view of the restrictions placed on the Stokes multipliers in \eqref{Stokes-multiplier-restriction}, we obtain
\begin{equation}\label{s2-s-1-asym}
s_{2}=-s_{-1}=-2A(\xi)e^{-2\xi E_{1}+2E_{2}+o(1)}\cos\left[2\xi E_{1}-2E_{3}+o(1)\right]
\end{equation}
and
\begin{equation}\label{s1-s0-s3-asym}
s_{1}=\overline{s_{3}}=-\overline{s_{0}}=\frac{e^{2\xi(1+i)E_{1}-2(E_{2}+i E_{3})+o(1)}}{A(\xi)}
\end{equation}
as $\xi\to+\infty.$

\begin{remark}\label{remark-case-I-negative}
If $A(\xi)\cdot B(\xi)<0$, the analysis is similar to the above argument, and the results are the same as the ones in \eqref{s2-s-1-asym} and \eqref{s1-s0-s3-asym}. To get the corresponding results, the scaling \eqref{scaling-caseI} should be replaced by $\lambda=i\xi^{\frac{1}{3}}\eta$. Moreover, one should also note the following fact:
\begin{equation}\label{a-integral}
\begin{split}
K:&=\frac{2A(\xi)^2-\left(B(\xi)/A(\xi)\right)^2}{8}\int_{\alpha_{1}}^{\infty}\left[\frac{1}{(s+\frac{B(\xi)}{2A(\xi)})}+\frac{1}{(s-\frac{B(\xi)}{2A(\xi)})}\right]\frac{1}{\sqrt{s^4+\frac{1}{4}}}ds\\
&=-\log|A(\xi)|+\frac{\pi i}{4},
\end{split}
\end{equation}
where $\alpha_{1}=\frac{1}{\sqrt{2}}e^{-\frac{\pi i}{4}}$ and $A(\xi)^4-B(\xi)^2=1$. The proof of \eqref{a-integral} is left in Appendix C.
\end{remark}

\subsection{Case two: $A(\xi)^4-B(\xi)^2=-1$}
\begin{figure}[h]
 \begin{center}
\setlength{\unitlength}{1cm}
\begin{tikzpicture}
%(4.3, 3.6)(-2.5,-0.25)
\qbezier(1,0)(1.314,1.00)(2,1.4)[thick];
\qbezier(1,0)(1.314,-1.00)(2,-1.4)[thick];
\draw(1,0)--(-1,0)[thick];
\qbezier(0,1)(0.768,0.568)(2,1.45)[thick];
\draw(0,1)--(0,2.5)[thick];
\qbezier(0,1)(-0.768,0.568)(-2,1.45)[thick];
\qbezier(-1,0)(-1.314,1.00)(-2,1.4)[thick];
\qbezier(-1,0)(-1.314,-1.00)(-2,-1.4)[thick];
\qbezier(0,-1)(0.768,-0.568)(2,-1.45)[thick];
\draw(0,-1)--(0,-2.5)[thick];
\qbezier(0,-1)(-0.768,-0.568)(-2,-1.45)[thick];
\linethickness{1.175mm}
\put( 1.1,0){$\hat{\eta}_1$};
\put( 0.12,1.1){$\hat{\eta}_2$};
\put( -1.45,0){$\hat{\eta}_3$};
\put(0.12,-1.3){$\hat{\eta}_4$};
\put(2.1,-1.5){$\gamma_{-1}$};
\put(2.1,1.5){$\gamma_{0}$};
\put(0,2.6){$\gamma_{1}$};
\put(-2.5,-1.5){$\gamma_{3}$};
\put(-2.5,1.5){$\gamma_{2}$};
\put(0,-2.7){$\gamma_{-2}$};
\put(0,0){\circle*{0.01}};
\put(1,0){\circle*{0.1}}[color=blue];
\put(0,1){\circle*{0.1}}[color=blue];
\put(-1,0){\circle*{0.1}}[color=blue];
\put(0,-1){\circle*{0.1}}[color=blue];
\put(-0.20,-0.4){$O$};
\end{tikzpicture}
\end{center}
 \caption{\small{\,Stokes geometry of $\hat{G}(\eta,\xi)d\eta^2$\,}.}
 \label{figure-case-II}
\end{figure}
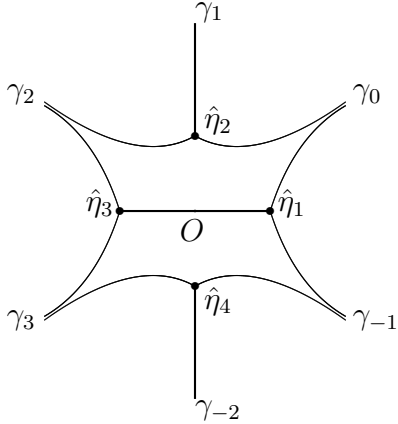

Scale $\lambda=-i\xi^{\frac{1}{3}}\eta$ and set $Y(\eta,\xi)=\phi(\lambda,\xi)$. In this case, the second order equation \eqref{Sec order equ2} turns to
\begin{equation}
\begin{split}\label{Second order equation asy 2}
  \frac{d^2Y}{d\eta^2}&=\xi^2\left\{16(\eta^4-\frac14)-\frac{8\eta}{\xi}+\frac{1}{\xi}\frac{8A(\xi)\eta^2-4A(\xi)^3}{2A(\xi)\eta+B(\xi)}   +\hat{g}(\eta,\xi)\right\}Y(\eta,\xi)\\
  &:=\xi^2\hat{G}(\eta,\xi)Y(\eta,\xi),
\end{split}
\end{equation}
where $\hat{g}(\eta,\xi)=\mathcal{O}\left(\xi^{-\frac{4}{3}}\right)$ as $\xi\to+\infty$ uniformly for all $\eta$ provided that $2A(\xi)\eta+B(\xi)$ is bounded away from $0$.  Note that there are four simple turning points, say $\hat{\eta}_j$($j=1,2,3,4$), near $\hat{\alpha}_1=\frac{1}{\sqrt{2}}$, $\hat{\alpha}_2=\frac{1}{\sqrt{2}}e^{\frac{\pi i}{2}}$, $\hat{\alpha}_3=-\frac{1}{\sqrt{2}}$ and $\hat{\alpha}_4=\frac{1}{\sqrt{2}}e^{-\frac{\pi i}{2}}$ respectively. It follows from \eqref{Second order equation asy 2} that $\frac{1}{\hat{\eta}_j-\hat{\alpha}_j}=\mathcal{O}(\xi)~(j=1,2,3,4)$ as $\xi\to+\infty$.
According to \cite{Eremenko}, the limiting state of the Stokes geometry of the quadratic form $\hat{G}(\eta,\xi)d\eta^2$ as $\xi\to+\infty$
is described in Figure \ref{figure-case-II}. Therefore, following \cite[Theorem 2]{APC}, we can approximate
the solutions of \eqref{Second order equation asy 2} via Airy functions.

Define two conformal mappings $\rho(\eta)$ and $\tau(\eta)$ by
\begin{equation}\label{Define rho}
\int_0^\rho s^{1/2}ds=\int_{\hat{\eta}_0}^{\eta}\hat{G}(s,\xi)^{1/2}d s,
\end{equation}
and
\begin{equation}\label{Define tau}
  \int_0^{\tau} s^{1/2}ds=\int_{\hat{\eta}_2}^\eta \hat{G}(s,\xi)^{1/2}ds,
\end{equation}
respectively near the neighborhoods of $\eta=\hat{\eta}_1$ and $\eta=\hat{\eta}_2$. Then the following lemma is a consequence of \cite[Theorem 2]{APC}.
%%%%%%%%%%%%%%%%%%%%%%%%%%%%%%%%%%%%%%%%%%%%%%%%%%%%%%%%%%%%%%%%%%%%%%%%%%%%%%%%%%%%%%%%%%%%%%%%%%%%%%%%%%%%%%%%%%%%%%%%%%

\begin{lemma}\label{Uniform asymptotic theorem2}
Given any solution $\phi(\eta,\xi)$ of \eqref{Second order equation asy 2},
there exist constants $d_1,~d_2$ and $\tilde{d}_1,~\tilde{d}_2$ such that
 \begin{equation}\label{Uniform appr phi21}
   \phi(\eta,\xi)=\Big(\frac{\rho}{\hat{G}(\eta,\xi)}\Big)^{\frac{1}{4}}
\Big\{[d_1+o(1)]{\rm Ai}(\xi^{\frac{2}{3}}\rho)
+[d_2+o(1)]{\rm Bi}(\xi^{\frac{2}{3}}\rho)\Big\},~\xi\to +\infty
 \end{equation}
 uniformly for $\eta$ on any two adjacent Stokes lines emanating from $\hat{\eta}_0$; and
\begin{equation}\label{Uniform appr phi22}
   \phi(\eta,\xi)=\Big(\frac{\tau}{\hat{G}(\eta,\xi)}\Big)^{\frac{1}{4}}
\Big\{[\tilde{d}_1+o(1)]{\rm Ai}(\xi^{\frac{2}{3}}\tau)
+[\tilde{d}_2+o(1)]{\rm Bi}(\xi^{\frac{2}{3}}\tau)\Big\},~\xi\to +\infty
 \end{equation}
 uniformly for $\eta$ on any two adjacent Stokes lines emanating from $\hat{\eta}_2$.
\end{lemma}
%%%%%%%%%%%%%%%%%%%%%%%%%%%%%%%%%%%%%%%%%%%%%%%%%%%%%%%%%%%%%%%%%%%%%%%%%%%%%%%%%%%%%%%%%%%%%%%%%%%%%%%%%%%%%%%%%
To calculate the connection martix $S_1$ as $\xi\rightarrow+\infty$, it is necessary to clarify the asymptotic behavior of
 $\rho(\eta)$ and $\tau(\eta)$ as $\xi,\eta\to\infty$.
%%%%%%%%%%%%%%%%%%%%%%%%%%%%%%%%%%%%%%%%%%%%%%%%%%%%%%%%%%%%%%%%%%%%
\begin{lemma}\label{lem relation rho tau and eta}
For large $\xi$ and $\eta$ with $|\eta|\gg\xi$,
\begin{equation}\label{relation-rho-eta-infty}
  \frac{2}{3}\rho^{\frac{3}{2}}=\frac{4}{3}\eta^3-\sqrt{2}E_{1} -\frac{1}{\xi}\ln(2\eta)+\frac{1}{2\xi}\log\frac{2A(\xi)\eta+B(\xi)}{B(\xi)}+\frac{F_{1}}{\xi}+o(\xi^{-1}),
  \end{equation}
and
\begin{equation}\label{relation-tau-eta-infty}
  \frac{2}{3}\tau^{\frac{3}{2}}=\frac{4}{3}\eta^3+\sqrt{2}iE_{1}-\frac{1}{\xi}\ln(2\eta)+\frac{1}{2\xi}\log(2A(\xi)\eta+B(\xi))+\frac{\pi i}{2\xi}+\frac{iF_{2}}{\xi}+o(\xi^{-1}),
\end{equation}
where $E_{1}$ is the same as in Lemma \ref{lem relation zeta and eta infty} and
\begin{equation}\label{def-F1}
F_{1}=\frac{1}{2}\int_{\hat{\alpha}_{1}}^{\infty}\frac{A(\xi)}{B(\xi)}\cdot \frac{1-A(\xi)^2}{2\frac{A(\xi)}{B(\xi)}s+1}\cdot \frac{1}{\sqrt{s^4-\frac{1}{4}}}ds-\frac{1}{2}\int_{\hat{\alpha}_{1}}^{\infty}\frac{s\log(2\frac{A(\xi)}{B(\xi)}s+1)}{\left(s^2+\frac{1}{2}\right)\sqrt{s^4-\frac{1}{4}}}ds,
\end{equation}
and
\begin{equation}\label{def-F2}
F_{2}=\frac{i}{2}\int_{\hat{\alpha}_{2}}^{\infty} \frac{A(\xi)}{B(\xi)}\cdot\frac{1+A(\xi)^2}{2\frac{A(\xi)}{B(\xi)}s+1}\cdot\frac{1}{\sqrt{s^4-\frac{1}{4}}}ds-\frac{i}{2}\int_{\hat{\alpha}_{2}}^{\infty}\frac{s\log(2\frac{A(\xi)}{B(\xi)}s+1)}{\left(s^2-\frac{1}{2}\right)\sqrt{s^4-\frac{1}{4}}}ds+\frac{i\log{B(\xi)}}{2}.
\end{equation}
Specially, if $A(\xi)=0$, then $F_{1}=0$, and $F_{2}=0(F_{2}=-\frac{\pi}{2})$ when $B(\xi)=1(B(\xi)=-1)$.
\end{lemma}
\begin{remark}
A more careful calculation to \eqref{def-F2} yields $\im F_{2}=0$ when $A(\xi)^4-B(\xi)^2=-1$. Actually, we see from \eqref{def-F2} that
\begin{equation}
\begin{split}
\im F_{2}&=-\int_{1}^{+\infty}\frac{(A(\xi)^2+A(\xi)^4)}{(2A(\xi)^2s+B(\xi)^2)\sqrt{s^2-1}}ds+\frac{1}{2}\int_{1}^{+\infty}\frac{\ln(2A(\xi)^2s+B(\xi)^2)}{(s+1)\sqrt{s^2-1}}ds\\
&=-\frac{A(\xi)^2}{A(\xi)^2-1}\log{A(\xi)^2}+\frac{A(\xi)^2}{A(\xi)^2-1}\log{A(\xi)^2}=0.
\end{split}
\end{equation}
\end{remark}

The proof of Lemma \ref{lem relation rho tau and eta} is left in Appendix B. Now we turn to prove Lemma
\ref{stokes-multipliers-case-II}, {\it i.e.} to calculate the Stokes multipliers $s_{-1}$ and $s_{0}$ as $\xi\rightarrow+\infty$ with $A(\xi)^4-B(\xi)^2=-1$.

\textbf{Proof of Lemma \ref{stokes-multipliers-case-II}.} It can be seen from Fig.\ref{figure-case-II} that, in order to calculate $s_{-1}$, we should consider the uniform asymptotic behaviors of  $\Ai(\xi^{\frac{2}{3}}\rho)$ and $\Bi(\xi^{\frac{2}{3}}\rho)$ on the two adjacent Stokes lines $\gamma_{-1}$ and $\gamma_{0}$. When $\arg\lambda\sim-\frac{2\pi}{3}$, we have $\arg\eta=\arg(i\lambda)\sim-\frac{\pi}{6}$, which implies $\arg\rho\sim -\frac{\pi}{3}$ as $\eta\to\infty$.
Applying the asymptotic behavior of $\Ai(z)$ and $\Bi(z)$ in \eqref{eq-behavior-Ai} and \eqref{eq-behavior-Bi},  we get
\begin{equation}\label{eq-Phi--pi/6}
\begin{cases}
\sqrt{2\lambda a-ib}\left(\frac{\rho}{\hat{G}(\eta,\xi)}\right)^{\frac{1}{4}}\Ai(\xi^{\frac{2}{3}}\rho)&\sim d_{1}e^{\frac{4i}{3}\lambda^{3}},\\
\sqrt{2\lambda a-ib}\left(\frac{\rho}{\hat{G}(\eta,\xi)}\right)^{\frac{1}{4}}\Bi(\xi^{\frac{2}{3}}\rho)&\sim -id_{1}e^{\frac{4i}{3}\lambda^{3}}+2d_{2}\frac{a}{2\lambda}e^{-\frac{4i}{3}\lambda^{3}},
\end{cases}
\end{equation}
where
\begin{equation}\label{eq-d1-d2}
d_{1}=\frac{\xi^{\frac{1}{6}}B(\xi)^{\frac{1}{2}}e^{-\frac{\pi i}{4}+\xi\sqrt{2}E_{1}-F_{1}+o(1)}}{2\sqrt{\pi}},\qquad d_{2}=-\frac{\xi^{\frac{1}{6}}B(\xi)^{-\frac{1}{2}}e^{\frac{\pi i}{4}-\xi\sqrt{2}E_{1}+F_{1}+o(1)}}{2\sqrt{\pi}}.
\end{equation}
as $\xi\rightarrow+\infty$. Comparing \eqref{eq-Phi--pi/6} to \eqref{eq-asy-Phi21-Phi22}, we have
\begin{equation}\label{Phi21-Phi22-to-Airy-gamma-1-2}
((\Phi_{-1})_{21},(\Phi_{-1})_{22})=\sqrt{2\lambda a-ib}\left(\frac{\rho}{\hat{G}(\eta,\xi)}\right)^{\frac{1}{4}}(\Ai(\xi^{\frac{2}{3}}\rho),\Bi(\xi^{\frac{2}{3}}\rho))
\left(\begin{matrix}
\frac{i}{2d_{2}}&\frac{1}{d_{1}}\\[0.2cm]
\frac{1}{2d_{2}}&0
\end{matrix}\right)
\end{equation}
as $\xi\to+\infty$. Here and after, the $d_{j}$'s in (\ref{Phi21-Phi22-to-Airy-gamma0}) are also not equal but asymptotically equal to the corresponding ones in (\ref{eq-d1-d2}) as $\xi\rightarrow+\infty$. By abuse of notations, we use the same symbol for the $d_{j}$'s in these two formulas, since we only care about the asymptotic behavior of the Stokes multipliers.

Similarly, when $\arg\lambda\sim-\frac{\pi}{3}$, we have $\arg\eta=\arg(i\lambda)\sim\frac{\pi}{6}$, which implies $\arg\rho\sim\frac{\pi}{3}$. Making use of the uniform asymptotic behaviors of $\Ai(z)$ and $\Bi(z)$ again and noting \eqref{eq-asy-Phi21-Phi22}, one can immediately obtain
\begin{equation}\label{Phi21-Phi22-to-Airy-gamma0-2}
\left((\Phi_0)_{21},(\Phi_0)_{22}\right)=\sqrt{2\lambda a-ib}\left(\frac{\rho}{\hat{G}(\eta,\xi)}\right)^{\frac{1}{4}}\left(\Ai(\xi^{\frac{2}{3}}\rho),\Bi(\xi^{\frac{2}{3}}\rho)\right)
\left(\begin{matrix}-\frac{i}{2d_{2}} & \frac{1}{d_{1}} \\ \frac{1}{2d_{2}}& 0\end{matrix}\right)
\end{equation}
as $\xi\rightarrow+\infty$. Hence,
\begin{equation}
\begin{split}
\left((\Phi_0)_{21},(\Phi_0)_{22}\right)&=((\Phi_{-1})_{21},(\Phi_{-1})_{22})\left(\begin{matrix}
\frac{i}{2d_{2}}&\frac{1}{d_{1}}\\[0.2cm]
\frac{1}{2d_{2}}&0
\end{matrix}\right)^{-1}\left(\begin{matrix}-\frac{i}{2d_{2}} & \frac{1}{d_{1}} \\ \frac{1}{2d_{2}}& 0\end{matrix}\right)\\
&=((\Phi_{-1})_{21},(\Phi_{-1})_{22})\left(\begin{matrix}
1&0\\[0.2cm]
-i\frac{d_{1}}{d_{2}}&0
\end{matrix}\right).
\end{split}
\end{equation}
This implies $s_{-1}=-i\frac{d_{1}}{d_{2}}=-B(\xi)e^{2\sqrt{2}E_{1}\xi-2F_{1}+o(1)}$ as $\xi\to+\infty$.

To derive $s_{0}$, we need to analyze the uniform asymptotic behaviors of $\Ai(\xi^{\frac{2}{3}}\tau)$ and $\Bi(\xi^{\frac{2}{3}}\tau)$ on the two adjacent Stokes lines $\gamma_{0}$ and $\gamma_{1}$ as $\xi,|\eta|\rightarrow+\infty$. When $\eta\to\infty$ on $\gamma_{0}$, we find that $\arg\eta=\frac{\pi}{6}$, which implies $\arg\tau\sim\frac{\pi}{3}$. Then according to \eqref{eq-behavior-Ai} and \eqref{eq-behavior-Bi}, we obtain
\begin{equation}\label{eq-Phi-pi/6-2}
\begin{cases}
\sqrt{2\lambda a-ib}\left(\frac{\tau}{\hat{G}(\eta,\xi)}\right)^{\frac{1}{4}}\Ai(\xi^{\frac{2}{3}}\tau)&\sim \tilde{d}_{1}e^{\frac{4i}{3}\lambda^{3}},\\
\sqrt{2\lambda a-ib}\left(\frac{\tau}{\hat{G}(\eta,\xi)}\right)^{\frac{1}{4}}\Bi(\xi^{\frac{2}{3}}\tau)&\sim i\tilde{d}_{1}e^{\frac{4i}{3}\lambda^{3}}+2\tilde{d}_{2}\frac{a}{2\lambda}e^{-\frac{4i}{3}\lambda^{3}},
\end{cases}
\end{equation}
where
\begin{equation}\label{eq-c1-c2}
\tilde{d}_{1}=\frac{\xi^{\frac{1}{6}}e^{-\frac{\pi i}{4}-\sqrt{2}i\xi E_{1}-\frac{\pi i}{2}-iF_{2}+o(1)}}{2\sqrt{\pi}},\qquad \tilde{d}_{2}= \frac{\xi^{\frac{1}{6}}e^{\frac{\pi i}{4}+\sqrt{2}i\xi E_{1}+\frac{\pi i}{2}+iF_{2}+o(1)}}{2\sqrt{\pi}}.
\end{equation}
as $\xi\rightarrow+\infty$. It can be further derived by matching \eqref{eq-Phi-pi/6-2} to \eqref{eq-asy-Phi21-Phi22} that
\begin{equation}\label{Phi21-Phi22-to-Airy-gamma0-3}
\left((\Phi_0)_{21},(\Phi_0)_{22}\right)=\sqrt{2\lambda a-ib}\left(\frac{\tau}{\hat{G}(\eta,\xi)}\right)^{\frac{1}{4}}\left(\Ai(\xi^{\frac{2}{3}}\tau),\Bi(\xi^{\frac{2}{3}}\tau)\right)
\left(\begin{matrix}-\frac{i}{2\tilde{d}_{2}} & \frac{1}{\tilde{d}_{1}} \\ \frac{1}{2\tilde{d}_{2}}& 0\end{matrix}\right)
\end{equation}
as $\xi\rightarrow+\infty$. When $\eta\rightarrow\infty$ on $\gamma_{1}$, we know that $\arg\eta\sim\frac{\pi}{2}$, then $\arg\tau\sim\pi$. In a similar way to the above argument, we have
\begin{equation}\label{Phi21-Phi22-to-Airy-gamma1-3}
\left((\Phi_1)_{21},(\Phi_1)_{22}\right)=\sqrt{2\lambda a-ib}\left(\frac{\tau}{\hat{G}(\eta,\xi)}\right)^{\frac{1}{4}}\left(\Ai(\xi^{\frac{2}{3}}\tau),\Bi(\xi^{\frac{2}{3}}\tau)\right)
\left(\begin{matrix}-\frac{i}{2\tilde{d}_{2}} & \frac{1}{2\tilde{d}_{1}} \\ \frac{1}{2\tilde{d}_{2}}& - \frac{i}{2\tilde{d}_{1}}\end{matrix}\right).
\end{equation}
A combination of \eqref{Phi21-Phi22-to-Airy-gamma0-3} and \eqref{Phi21-Phi22-to-Airy-gamma1-3} yields
\begin{equation}\label{Stokes-matrix-s0-II}
\begin{split}
\left((\Phi_1)_{21},(\Phi_1)_{22}\right)&=\left((\Phi_0)_{21},(\Phi_0)_{22}\right)\left(\begin{matrix}-\frac{i}{2\tilde{d}_{2}} & \frac{1}{\tilde{d}_{1}} \\ \frac{1}{2\tilde{d}_{2}}& 0\end{matrix}\right)^{-1}
\left(\begin{matrix}-\frac{i}{2\tilde{d}_{2}} & \frac{1}{2\tilde{d}_{1}} \\ \frac{1}{2\tilde{d}_{2}}& - \frac{i}{2\tilde{d}_{1}}\end{matrix}\right)\\
&=\left((\Phi_0)_{21},(\Phi_0)_{22}\right)\left(\begin{matrix}1& -i\frac{\tilde{d}_{2}}{\tilde{d}_{1}} \\ 0& 1\end{matrix}\right).
\end{split}
\end{equation}
This means that $s_{0}=-i\frac{\tilde{d}_{2}}{\tilde{d}_{1}}=-e^{2\sqrt{2}\xi i E_{1}+2iF_{2}+o(1)}$ as $\xi\to+\infty$.
Finally, in view of the restriction of the Stokes multipliers in \eqref{Stokes-multiplier-restriction}, we obtain
\begin{equation}\label{s2-s-1-asym-II}
s_{2}=-s_{-1}=B(\xi)e^{2\sqrt{2}E_{1}\xi-2F_{1}+o(1)}
\end{equation}
and
\begin{equation}\label{s1-s0-s3-asym-II}
s_{1}=\overline{s_{3}}=-\overline{s_{0}}=e^{-2\sqrt{2}\xi i E_{1}-2iF_{2}+o(1)}
\end{equation}
as $\xi\to+\infty$.

\begin{remark}
When $A(\xi)\cdot B(\xi)<0$, the analysis is similar to the above argument, and the corresponding results are agreement with \eqref{s2-s-1-asym-II} and \eqref{s1-s0-s3-asym-II}.
\end{remark}

\section*{Acknowledgements}
The authors are grateful to the referees for their valuable suggestions and comments for improvement of the paper.
The work of Wen-Gao Long was supported by the Natural Science Foundation of Hunan Province under Grant Number 2020JJ5152, the General Project of Hunan Provincial Department of Education under Grant Number 19C0771 and the Doctoral Startup Fund of Hunan University of Science and Technology under Grant Number E51871. The work of Zhao-Yun Zeng was supported in part by the National Natural
Science Foundation of China under Grant Numbers 11571375 and 11571376, the Scientific and Technological Research Project of Jiangxi Provincial Education Department (numbers GJJ180581 and GJJ190549), Doctoral Startup Fund of Jinggangshan University (number JZB16001), and Open Research Fund Program of the State Key Laboratory of Low-Dimensional Quantum Physics (number KF201903).
%%%%%%%%%%%%%%%%%%%%%%%%%%%%%%%%%%%%%%%%%%%%%%%%%%%%%%%%%%%%%%%%%%%%%%%%%%%%%%%%%%%%%%%%%%%%%%%%%%%%%%%%
\section*{Appendix }
\begin{appendix}{}
 % \renewcommand{\thesection}{\Alph{section}}
%\appendix{}
%  \renewcommand{\appendixname}{Appendix~\Alph{section}}
%\section{Proof of Lemma \ref{lem relation zeta and eta infty} } \renewcommand{\thesection}{\Alph{section}}
\section{Proof of Lemma \ref{lem relation zeta and eta infty}}
The idea to prove Lemma \ref{lem relation zeta and eta infty} is to compute the asymptotic behavior of the integral on the right hand side of \eqref{Define zeta} as $\xi,\eta\to\infty$. Here we only show the validity of \eqref{eq-vartheta(eta)-infty}.
The proof of \eqref{eq-zeta(eta)-infty} is similar and hence omitted here.

Choose $T^*$ to be a point such that $|T^*-\eta_1|\sim \xi^{-\frac 3 4}$, then
 \begin{equation}\label{split right integral1}
   \int_{\eta_1}^\eta G(s,\xi)^{\frac 1 2}ds=\int_{\eta_1}^{T^*} G(s,\xi)^{\frac 1 2}ds
   +\int_{T^*}^\eta G(s,\xi)^{\frac 1 2}ds:=I_1+I_2.
 \end{equation}
Since $|T^*-\eta_1|\sim \xi^{-\frac 3 4}$ and
$\frac{1}{\eta_1-\alpha_1}=\mathcal{O}\left(\xi\right)$ as $\xi\to+\infty$, then $|T^*-\alpha_1|\sim\xi^{-\frac 3 4}$ uniformly for all $s$ on the integral contour of $I_{2}$. It implies that $I_1=o(\xi^{-1})$ as $\xi\to\infty$. By the Taylor
expansion of $G(s,\xi)$ with respect to $\frac{1}{\xi}$, we have
\begin{equation}\label{expand F1/2}
   G(s,\xi)^{\frac 1 2}=G_{0}(s,\xi)-\frac{1}{\xi}G_{1}(s,\xi)+\mathcal{O}\left(\frac{1}{\xi^{\frac{4}{3}}\sqrt{s^4+\frac 1 4}}\right)
\end{equation}
 as $\xi\to+\infty$, where $G_{0}(s,\xi)=4\sqrt{s^4+\frac{1}{4}}$ and
\begin{equation}\label{eq-F1(s,xi)}
G_{1}(s,\xi)=\frac{s}{2\sqrt{s^4+\frac 1 4}}+\frac{B(\xi)}{4A(\xi)}\frac{1}{\sqrt{s^4+\frac{1}{4}}}+\frac{1}{8}\frac{2 A(\xi)^2-\left(B(\xi)/A(\xi)\right)^2}{(s+\frac{B(\xi)}{2 A(\xi)})\sqrt{s^4+\frac 1 4}}.
\end{equation}
Hence
\begin{equation}\label{eq-I2}
\begin{split}
I_{2}=\int_{T^*}^{\eta}\left[G_{0}(s,\xi)-\frac{1}{\xi}G_{1}(s,\xi)\right]ds+o(\xi^{-1})=\int_{\alpha_{1}}^{\eta}\left[G_{0}(s,\xi)-\frac{1}{\xi}G_{1}(s,\xi)\right]ds+o(\xi^{-1})
\end{split}\end{equation}
as $\xi\rightarrow+\infty$. To get the final approximation, one may only need to note the facts $|T^*-\eta_1|\sim \xi^{-\frac 3 4}$  and
$\frac{1}{\eta_1-\alpha_1}=\mathcal{O}\left(\xi\right)$ as $\xi\to+\infty$. Furthermore, integrating by parts, we have
\begin{align}\label{first of I2}
\int_{\alpha_1}^{\eta} G_{0}(s,\xi)ds&=\int_{\alpha_1}^\eta4\sqrt{s^4+\frac 1 4}ds=2\int_{\alpha_1}^\eta\sqrt{4s^4+1}ds\\
  &=2\eta\sqrt{4\eta^4+1}-2\int_{\alpha_1}^\eta4\sqrt{s^4+\frac 1 4}ds+4\int_{\alpha_1}^\eta\frac{1}{\sqrt{4s^4+1}}ds.
\end{align}
Hence, we get
\begin{equation}\label{approx-int-F0}
\begin{split}
\int_{\alpha_1}^{\eta} G_{0}(s,\xi)ds=&\frac{2}{3}\eta\sqrt{4\eta^4+1}+\frac{4}{3}\int_{\alpha_1}^\eta\frac{1}{\sqrt{4s^4+1}}ds\\
=&\frac{4}{3}\eta^{3}+\frac{4}{3}\int_{\alpha_1}^\infty\frac{1}{\sqrt{4s^4+1}}ds+\mathcal{O}\left(\frac{1}{\eta}\right)\\
=&\frac{4}{3}\eta^{3}+(1+i)E_{1}+o(\xi^{-1})
\end{split}
\end{equation}
as $\xi,\eta\to\infty$ provided that $|\eta|\gg \xi$, where
\begin{equation}\label{eq-E1}
E_{1}:=\frac{4}{3(1+i)}\int_{\alpha_1}^\infty\frac{1}{\sqrt{4s^4+1}}ds=\frac{1}{6}B\left(\frac{1}{4},\frac{1}{2}\right).
\end{equation}
To derive the asymptotic behavior of $\int_{\alpha_{1}}^{\eta}G_{1}(s,\xi)ds$ as $\xi\rightarrow+\infty$, we still need to calculate $\int_{\alpha_{1}}^{\eta}\frac{s}{\sqrt{s^4+1/4}}ds$. As a matter of fact,
\begin{equation}\label{eq-1/(s4+1)}
\int_{\alpha_{1}}^{\eta}\frac{s}{\sqrt{s^4+\frac{1}{4}}}ds=\frac{1}{2}\ln(s^2+\sqrt{s^4+\frac 1 4})\Big|_{\alpha_{1}}^{\eta}=\ln(2\eta)+\frac{\pi i}{4}+o(1)
\end{equation}
as $\xi,\eta\to+\infty$ with $\eta\gg\xi$. Combining \eqref{eq-F1(s,xi)}, \eqref{eq-E1} and \eqref{eq-1/(s4+1)}, we have
\begin{equation}\label{eq-F1-behavior}
\int_{\alpha_{1}}^{\eta}G_{1}(s,\xi)ds=\frac{1}{2}\ln(2\eta)+E_{2}+iE_{3}+o(1)
\end{equation}
as $\xi,\eta\to+\infty$ with $\eta\gg\xi$ and $\arg\eta\in\left(-\frac{\pi}{2},\frac{\pi}{6}\right)$, where $E_{2}=\frac{B(\xi)}{16A(\xi)}B\left(\frac{1}{4},\frac{1}{2}\right)+\Re\mathcal{F}(\xi)$,  $E_{3}=\frac{B(\xi)}{16A(\xi)}B\left(\frac{1}{4},\frac{1}{2}\right)+\frac{\pi}{8}+\im\mathcal{F}(\xi)$, and
$$\mathcal{F}(\xi):=\frac{2 A(\xi)^2-\left(B(\xi)/A(\xi)\right)^2}{8}\int_{\alpha_{1}}^{\infty}\frac{1}{(s+\frac{B(\xi)}{2 A(\xi)})\sqrt{s^4+\frac 1 4}}ds.$$
A combination of \eqref{eq-I2}, \eqref{approx-int-F0} and \eqref{eq-F1-behavior} yields \eqref{eq-vartheta(eta)-infty}. Finally, noting the symmetry property $\vartheta(\eta)=\overline{\zeta(\bar{\eta})}+o(\xi^{-1})$ as $\xi,\eta\to\infty$, we get \eqref{eq-zeta(eta)-infty}.

\section{Proof of Lemma \ref{lem relation rho tau and eta}}\renewcommand{\thesection}{\Alph{section}}
The idea to prove Lemma \ref{lem relation rho tau and eta} is to compute the asymptotic behavior of the integral on the right hand side of \eqref{Define rho} as $\xi,\eta\to\infty$. Here we only show the validity of \eqref{relation-rho-eta-infty}.
The proof of \eqref{relation-tau-eta-infty} is similar and hence omitted here.

 Choose $T_1^*$ to be a point such that $|T_1^*-\hat{\eta}_0|\sim \xi^{-\frac 3 4}$, then
 \begin{equation}\label{split right integral2}
   \int_{\hat{\eta}_0}^\eta \hat{G}(s,\xi)^{\frac 1 2}ds=\int_{\hat{\eta}_0}^{T_1^*} \hat{G}(s,\xi)^{\frac 1 2}ds
   +\int_{T_1^*}^\eta \hat{G}(s,\xi)^{\frac 1 2}ds:=J_1+J_2.
 \end{equation}
Note the fact that $\hat{G}(s,\xi)=\mathcal{O}(\xi^{\frac{3}{4}})$ as $\xi\to+\infty$ uniformly on the integration contour of $J_1$. Then an easy approximation gives that $J_1=o(\xi^{-1})$ as $\xi\to\infty$. Since $|T_1^*-\hat{\eta}_1|\sim \xi^{-\frac 3 4}$ and
$\frac{1}{\hat{\eta}_1-\hat{\alpha}_1}=\mathcal{O}(\xi)$ as $\xi\to+\infty$, then $|T_1^*-\hat{\alpha}_1|>c\xi^{-\frac 3 4}$
for some constant $c>0$. By the Taylor expansion of $\hat{G}(s,\xi)$ with respect to $\frac{1}{\xi}$, we have
 \begin{equation}\label{expand Q1/2}
   \hat{G}(s,\xi)^{\frac 1 2}=4\sqrt{s^4-\frac 1 4}-\frac{1}{\xi}\frac{s}{\sqrt{s^4-\frac 1 4}}+\frac{1}{2\xi}\frac{2A(\xi)s^2-A(\xi)^3}{2A(\xi)s+B(\xi)}\frac{1}{\sqrt{s^4-\frac{1}{4}}}
   +\mathcal{O}\left(\frac{1}{\xi^{\frac{4}{3}}\sqrt{s^4-\frac 1 4}}\right)
 \end{equation}
as $\xi\to+\infty$. Making use of the fact $|T_1^*-\hat{\alpha}_1|>c\xi^{-\frac 3 4}$ again and noting that $\frac{1}{\sqrt{s^4-1}}=\mathcal{O}\left((s-\hat{\alpha}_{1})^{-\frac{1}{2}}\right)$ as $s\to\hat{\alpha}_{1}$, we find that
\begin{equation}\label{J2-reduced}
\begin{split}
J_{2}=&\int_{\hat{\alpha}_{1}}^{\eta}4\sqrt{s^4-\frac 1 4}ds-\frac{1}{\xi}\int_{\hat{\alpha}_{1}}^{\eta}\frac{s}{\sqrt{s^4-\frac 1 4}}ds +\frac{1}{2\xi}\int_{\hat{\alpha}_{1}}^{\eta}\frac{2A(\xi)s^2-A(\xi)}{2A(\xi)s+B(\xi)}\frac{1}{\sqrt{s^4-\frac{1}{4}}}ds\\
&+\frac{1}{2\xi}\int_{\hat{\alpha}_{1}}^{\infty}\frac{A(\xi)-A(\xi)^3}{2A(\xi)s+B(\xi)}\frac{1}{\sqrt{s^4-\frac{1}{4}}}ds+o(\xi^{-1})
\end{split}
\end{equation}
as $\xi,\eta\to\infty$ with $|\eta|\gg \xi$.
The the first integral of the right hand side in \eqref{J2-reduced} can be calculated as follows:
\begin{align}\label{first of J21}
2\int_{\hat{\alpha}_1}^\eta\sqrt{4s^4-1}ds=2\eta\sqrt{4\eta^4-1}-4\int_{\hat{\alpha}_1}^\eta\sqrt{4s^4-1}ds-4\int_{\hat{\alpha}_1}^\eta\frac{1}{\sqrt{4s^4-1}}ds,
\end{align}
which gives that
\begin{equation}\label{expr hatJ}
\begin{split}
2\int_{\hat{\alpha}_1}^\eta\sqrt{4s^4-1}ds&=\frac{2}{3}\eta\sqrt{4\eta^4-1}-\frac{4}{3}\int_{\hat{\alpha}_1}^\eta\frac{1}{\sqrt{4s^4-1}}ds\\
&=\frac{4}{3}\eta^3-\frac{4}{3}\int_{\hat{\alpha}_1}^{\infty}\frac{1}{\sqrt{4s^4-1}}ds
\end{split}
\end{equation}
as $\xi,\eta\to\infty$ provided that $|\eta|\gg \xi$. Next, we compute the integral
 $\int_{\hat{\alpha}_1}^{\infty}\frac{1}{\sqrt{4s^4-1}}ds$ explicitly. Set $s=\frac{t}{\sqrt{2}}$, then we have
\begin{equation}\label{B euler intergal trans2}
\int_{\hat{\alpha}_1}^{\infty}\frac{1}{\sqrt{4s^4-1}}ds=\frac{1}{\sqrt{2}}\int_1^{\infty}\frac{1}{\sqrt{t^4-1}}dt
=\frac {1}{4\sqrt{2}} B(\frac 1 2,\frac 1 4).
\end{equation}
The last equality can be derived by setting $u=t^{-4}$ and calculate as
\begin{equation}\label{B euler int2}
  \int_1^{\infty}\frac{1}{\sqrt{t^4-1}}dt=\frac{1}{4}\int_0^1(1-u)^{-\frac12}u^{-\frac34}du
  =\frac14 B(\frac 12,\frac 14).
\end{equation}
Then, combining \eqref{expr hatJ} and \eqref{B euler intergal trans2}, we obtain
\begin{equation}\label{first int value J2}
  \int_{\hat{\alpha}_{1}}^\eta4\sqrt{s^4-\frac 1 4}ds=\frac{4}{3}\eta^3-\frac{\sqrt{2}}{6} B(\frac 1 2,\frac 1 4)+o(\xi^{-1})
\end{equation}
as $\xi,\eta\to\infty$ provided that $|\eta|\gg\xi$.
The second integral on the right hand side of \eqref{J2-reduced} is
\begin{align}\label{Second int value J2}
  \int_{\hat{\alpha}_{1}}^\eta\frac{s}{\sqrt{s^4-\frac 1 4}}ds&=\frac{1}{2}\int_{\hat{\alpha}_{1}}^\eta\frac{ds^2}{\sqrt{s^4-\frac 1 4}}
  =\frac{1}{2}\ln(s^2+\sqrt{s^4-\frac 1 4})\Big|_{\hat{\alpha}_{1}}^\eta=\ln(2\eta)+o(1).
\end{align}
By integration by parts, we find that
\begin{equation}\label{Third int val J2}
\begin{split}
\int_{\hat{\alpha}_{1}}^{\eta}\frac{2A(\xi)\left(s^2-\frac{1}{2}\right)}{(2A(\xi)s+B(\xi))\sqrt{s^4-\frac{1}{4}}}ds=&\log(2A(\xi)\eta+B(\xi))-\int_{\hat{\alpha}_{1}}^{\infty}\frac{s\log(2A(\xi)s+B(\xi))}{\left(s^2+\frac{1}{2}\right)\sqrt{s^4-\frac{1}{4}}}ds+o(\xi^{-1})
\end{split}
\end{equation}
as $\xi\to+\infty$ with $|\eta|\gg\xi$.

Finally, combining \eqref{J2-reduced},
  \eqref{first int value J2}, \eqref{Second int value J2} and \eqref{Third int val J2}, we obtain
$$J_2=\frac{4}{3}\eta^3-\frac{\sqrt{2}}{6} B(\frac 1 2,\frac 1 4)-\frac{1}{\xi}\ln(2\eta)+\frac{1}{2\xi}\log(2A(\xi)\eta+B(\xi))+\frac{1}{\xi}F_{1}+o(\xi^{-1})$$
as $\xi\to+\infty$ with $|\eta|\gg\xi$, where
\begin{equation}\label{def-F1}
F_{1}=\frac{1}{2}\int_{\hat{\alpha}_{1}}^{\infty}\frac{A(\xi)-A(\xi)^3}{2A(\xi)s+B(\xi)}\frac{1}{\sqrt{s^4-\frac{1}{4}}}ds
-\frac{1}{2}\int_{\hat{\alpha}_{1}}^{\infty}\frac{s\log(2A(\xi)s+B(\xi))}{\left(s^2+\frac{1}{2}\right)\sqrt{s^4-\frac{1}{4}}}ds.
\end{equation}
This completes the proof of Lemma \ref{lem relation rho tau and eta}.

\section{Proof of \eqref{a-integral}}
In fact, with the transformation $t=2s^2$, we have
\begin{equation}
\begin{split}
&\int_{\alpha_{1}}^{\infty}\left[\frac{1}{(s+\frac{B(\xi)}{2A(\xi)})}+\frac{1}{(s-\frac{B(\xi)}{2A(\xi)})}\right]\frac{1}{\sqrt{s^4+\frac{1}{4}}}ds\\
=&2\int_{e^{-\frac{\pi i}{2}}}^{\infty}\frac{dt}{\left(t-\left(\frac{B(\xi)}{\sqrt{2}A(\xi)}\right)^2\right)\sqrt{t^2+1}}\\
=&\frac{4A(\xi)^2}{\sqrt{B(\xi)^4+4A(\xi)^4}}\left(\frac{\pi i}{2}-2\log{|A(\xi)|}\right).
\end{split}
\end{equation}
Making use of the fact $A(\xi)^4-B(\xi)^2=1$, we further obtain
\begin{equation}
\begin{split}
K&=\frac{2A(\xi)^2-\left(B(\xi)/A(\xi)\right)^2}{8}\frac{4A(\xi)^2}{\sqrt{B(\xi)^4+4A(\xi)^4}}\left(\frac{\pi i}{2}-2\log{|A(\xi)|}\right)\\
&=\frac{\pi i}{4}-\log{|A(\xi)|}.
\end{split}
\end{equation}
\end{appendix}

%%%%%%%%%%%%%%%%%%%%%%%%%%%%%%%%%%%%%%%%%%%%%%%%%%%%%%%%%%%%%%%%%%%%%%%%%%%%%%%%%%%%%%%%%%

\end{document}